**Michiel Hazewinkel** 1 CWI
Direct line: +31-20-5924204 POBox 94079
Secretary: +31-20-5924233 1090GB Amsterdam
Fax: +31-20-5924166
E-mail: mich@cwi.nl

# Hopf algebras of endomorphisms of Hopf algebras

by
*Michiel Hazewinkel*
*CWI*
*POBox 94079*
*1090GB Amsterdam*
*The Netherlands*

**Abstract**. In the last decennia two generalizations of the Hopf algebra of symmetric functions have appeared and shown themselves important, the Hopf algebra of noncommutative symmetric functions *NSymm* and the Hopf algebra of quasisymmetric functions *QSymm*. It has also become clear that it is important to understand the noncommutative versions of such important structures as *Symm* the Hopf algebra of symmetric functions. Not least because the right noncommmutative versions are often more beautiful than the commutaive ones (not all cluttered up with counting coefficients). *NSymm* and *QSymm* are not truly the full noncommutative generalizations. One is maximally noncommutative but cocommutative, the other is maximally non cocommutative but commutative. There is a common, selfdual generalization, the Hopf algebra of permutations of Malvenuto, Poirier, and Reutenauer (*MPR*). This one is, I feel, best understood as a Hopf algebra of endomorphisms. In any case, this point of view suggests vast generalizations leading to the Hopf algebras of endomorphisms and word Hopf algebras with which this paper is concerned. This point of view also sheds light on the somewhat mysterious formulas of *MPR* and on the question where all the extra structure (such as autoduality) comes from. The paper concludes with a few sections on the structure of *MPR* and the question of algebra retractions of the natural inclusion of Hopf algebras $NSymm \longrightarrow MPR$ and coalgebra sections of the dual natural projection of Hopf algebras $MPR \longrightarrow QSymm$. Several of these will be described explicitly.

## 1. Introduction and motivation.

The original motivation for these studies comes from the *MPR* Hopf algebra (Hopf algebra of permutations) introduced and studied by Malvenuto, Poirier and Reutenauer, [10, 11]. Here is a partial very incomplete description. Much more will be said later in this paper.

As a graded Abelian group *MPR* has a basis consisting of the empty permutation and all permutations on $n$ letters, $n = 1, 2, \cdots$ . Thus

$$MPR = \mathbf{Z} \oplus \bigoplus_{n=1}^{\infty} \mathbf{Z}S_n \tag{1.1}$$

where $S_n$ is the symmetric group on $n$ letters. The underlying Abelain group of the group ring

$\mathbf{Z}S_n$ is the homogeneous summand of degree $n$. This includes the empty permutation which is the (canonical) basis element of the first summand of (1.1).

Permutations on $m$ letters are written as words, with a word $\alpha = [a_1, a_2, \cdots, a_m]$ on the alphabet $\{1, \cdots, m\}$ which has no repeats, corresponding to the permutation $i \mapsto a_i$, $i = 1, \cdots, m$. The empty permutation corresponds to the empty word $[\;]$.

*Multiplication on MPR.* With these notations the multiplication on $MPR$ can be described as follows. The empty word serves as the unit element, and if $\alpha = [a_1, \cdots, a_m]$ and $\beta = [b_1, \cdots, b_n]$ are two permutation words their product is

$$m(\alpha \otimes \beta) = [a_1, \cdots a_m] \times_{sh} [m + b_1, \cdots, m + b_n] \tag{1.2}$$

where $\times_{sh}$ stands for the shuffle product of two words, which can be described as follows. The shuffle product of two words $[c_1, \cdots, c_p]$ and $[d_1, \cdots, d_q]$ is the sum of all permutations of $[c_1, \cdots c_p, d_1, \cdots, d_q]$ (with multiplicities) for which all the $c$'s and all the $d$'s occur in their original order. For instance $[3] \times_{sh} [1,2,4] = [3,1,2,4] + [1,3,2,4] + [1,2,3,4] + [1,2,4,3]$ and $[1] \times_{sh} [1,2] = 2[1,1,2] + [1,2,1]$. Note that because of the 'shift by $m$' in (1.2) the right hand side of (1.2) is a sum of permutation words (and there are no multiplicities).

*Standardization.* To describe the comultiplication the notion of standardization is needed. For any word $\alpha = [a_1, a_2, \cdots, a_m]$ without repeats over the alphabet of natural numbers $\mathbf{N} = \{1,2,3, \cdots\}$ its standardization is the permutation word $\mathrm{st}(\alpha) = [\varphi(a_1), \cdots, \varphi(a_m)]$ where $\varphi : \{a_1, \cdots, a_m\} \longrightarrow \{1, \cdots, m\}$ is the unique order preserving map between these ordered sets. For instance $\mathrm{st}([5,2,1,8] = [3,2,1,4]$. This notion of standardization is a special case of standardization of words as introduced by Schensted, [13], which applies to all words over $\mathbf{N}$, not only words without repeats.

*The comultiplication on MPR.* The comultiplication on $MPR$ is now defined by

$$\mu(\alpha) = \sum_{\alpha' * \alpha'' = \alpha} st(\alpha') \otimes st(\alpha'') \tag{1.3}$$

where the sum is over all cuts of $\alpha$, that is all pairs of words $(\alpha', \alpha'')$ whose concatenation $\alpha' * \alpha''$ is equal to $\alpha$.

With $1 = [\;]$ as unit and a counit $\varepsilon$ defined by $\varepsilon(\mathbf{Z}S_n) = 0$ and $\varepsilon([\;]) = 1$, $MPR$ becomes a bialgebra; there is also (of course, given the connected graded setting) an antipode, making $MPR$ a Hopf algebra.

The question which intrigued me is now: How does one dream up such formulas and are there natural generalizations?

For instance, there is a very natural generalization of (1.2) to arbitrary words as follows. Define the height of a word $\alpha = [a_1, a_2, \cdots, a_m]$ as $\mathrm{ht}(\alpha) = \max\{a_1, \cdots, a_m\}$ and define the product of two arbitrary words $\alpha = [a_1, a_2, \cdots, a_m]$ and $\beta = [b_1, \cdots, b_n]$ as

$$m_{WHA}(\alpha \otimes \beta) = [a_1, \cdots, a_m] \times_{sh} [\mathrm{ht}(\alpha) + b_1, \cdots, \mathrm{ht}(\alpha) + b_n] \tag{1.4}$$

(Note that this agrees with (1.2) when $\alpha$ is a permutation word.) The question now arises whether there is a corresponding comultiplication turning the graded Abelian group with as basis all words over the natural numbers into a Hopf algebra (with (1.4) as its multiplication). As it turns out there is and the result is what I call the word Hopf algebra (*WHA*). Here are some examples of the comultiplication

$$\alpha = [3,2,7,2,4],\ \mu(\alpha) = 1\otimes\alpha + [1]\otimes[2,6,2,3] + [3,2,6,2]\otimes[1] + \alpha\otimes 1$$

$$\alpha = [7,3,2,2,4],$$
$$\mu(\alpha) = 1\otimes\alpha + [3]\otimes[3,2,2,4] + [4,1]\otimes[2,2,3] + [6,3,2,2]\otimes[1] + \alpha\otimes 1$$

and it does not seem easy (to me) to guess at these formulas, or more precisely at the general recipe behind them.

There are also many more generalizations such as *dWHA*, the double word Hopf algebra.

There are good reasons for looking for generalizations of *MPR*. It is a most important and elegant Hopf algebra. On the one hand it generalizes *NSymm*, the Hopf algebra of noncommutative symmetric functions in the sense that the latter is a sub Hopf algebra of *MPR*, on the other hand it generalizes *QSymm*, the Hopf algebra of quasisymmetric functions in the sense that the latter is a quotient Hopf algebra of *MPR*. Moreover *MPR* is selfdual (with respect to a nondegenerate but not positive definite, inner product) and this duality is compatible with the duality between *NSymm* and *QSymm*. On the other hand *MPR* seems a little small in the sense that there appears to be no room for Frobenius morphisms (such as on *QSymm*) or Verschiebung morphisms (such as on *NSymm*). On *WHA* and *dWHA* (and their many relatives) there is plenty room for such morphisms.

## 2. Some notations and conventions.

Unless otherwise stated all Hopf algebras, morphisms, constructions, ... below are over **Z**, the ring of integers. Unadorned tensor products are over the integers. The corresponding objects over other commutative rings $k$ with unit element are obtained by tensoring. E.g. $MPR_k = MPR\otimes_{\mathbf{Z}} k$. It will always be assumed that the underlying module of a Hopf algebra is free. Everything also works over arbitrary base rings (commutative with unit element).

A Hopf algebra $H = (H, m, \mu, e, \varepsilon, \iota)$ over $k$ is graded if there is a decomposition

$$H = \bigoplus_{n=0}^{\infty} H_n \tag{2.1}$$

such that the multiplication $m: H\otimes_k H \longrightarrow H$, the comultiplication $\mu: H \longrightarrow H\otimes_k H$, the unit morphism $e: k \longrightarrow H$, the counit morphism $\varepsilon: H \longrightarrow k$, and the antipode $\iota: H \longrightarrow H$ satisfy:

$$\begin{aligned} &m(H_n\otimes_k H_m) \subset H_{n+m},\quad \mu(H_n) \subset \bigoplus_{i+j=n} H_i\otimes_k H_j \\ &e(k)\subset H_0,\quad \varepsilon(H_n) = 0 \text{ for } n>0,\quad \iota(H_n)\subset H_n \end{aligned} \tag{2.2}$$

Unless something else is explicitly said the term graded Hopf algebra also is assumed to imply that the homogeneous summands $H_n$ of degree $n$ are free of finite rank, and that the graded Hopf algebra is connected, which means that $H_0$ is free of rank 1, so that $e$ and $\varepsilon$ induce isomorphisms (of trivial Hopf algebras) $k \cong H_0$.

A morphism of graded Hopf algebras $\varphi: H \longrightarrow K$ is homogeneous iff $\varphi(H_n)\subset K_n$.

If $H$ is a Hopf algebra of finite rank, $H^*$ denotes the dual Hopf algebra. For a graded Hopf algebra $H$, the graded dual is

$$H^{\mathrm{gr}\,*} = \bigoplus_{n=0}^{\infty} H_n^*,\quad H_n^* = \mathrm{Mod}_k(H_n, k) \tag{2.3}$$

(where $\mathrm{Mod}_k$ means taking module morphisms).

If $\alpha = [a_1, a_2, \cdots, a_m]$ is a word over the natural numbers $\mathbf{N} = \{1,2,3,\cdots\}$, its length is $\lg(\alpha) = m$, its weight is $\mathrm{wt}(\alpha) = \sum_{i=1}^{m} a_i$, its height is $\mathrm{ht}(\alpha) = \max\{a_1, \cdots, a_m\}$, its support is $\mathrm{supp}(\alpha) = \{a_1, \cdots, a_m\}$, i.e. the set of different letters occurring in $\alpha$, and its multisupport is the multiset $\mathrm{msupp}(\alpha) = \{a_1, \cdots, a_m\}$, i.e. the different letters of $\alpha$ together with their multiplicities.

**3. The Hopf algebra *LieHopf* and its graded dual *Shuffle*.**

The elements of the Hopf algebra *MPR* can be interpreted as endomorphisms (of the underlying Abelian group) of the Hopf algebra *Shuffle* and this is a good way of looking at them which permits (natural) generalizations. So, here is a brief description of *LieHopf* and its graded dual *Shuffle.*

As an algebra *LieHopf* is simply the free associative algebra over the integers in countably many (noncommuting) indeterminates $U_1, U_2, \cdots$

$$LieHopf = \mathbf{Z}\langle U_1, U_2, \cdots \rangle \tag{3.1}$$

The comultiplication is determined by

$$\mu(U_n) = 1 \otimes U_n + U_n \otimes 1 \tag{3.2}$$

and the requirement that $\mu$ be an algebra morphism. The counit takes each $U_n$ to zero (and 1 to 1) and the antipode is determined by $\iota(U_n) = -U_n$ (and the requirement that it has to be an antimorphism of algebras).

*LieHopf* is graded with the grading determined by $\deg(U_n) = n$.

Thus a basis of *LieHopf* (as an Abelian group) is formed by all words in the alphabet $\{U_1, U_2, \cdots,\}$ and the multiplication is simply concatenation. If $\alpha = [a_1, a_2, \cdots, a_m]$ is a word over the natural numbers $U_\alpha$ is short for

$$U_\alpha = U_{a_1} U_{a_2} \cdots U_{a_m}, \quad U_{[\,]} = 1 \tag{3.3}$$

The degree of a noncommutative monomial $U_\alpha$ is the weight of the word $\alpha$.

A subword of $\alpha$ is any word of the form $\beta = [a_{i_1}, a_{i_2}, \cdots, a_{i_n}]$ with $i_1 < i_2 < \cdots < i_n$. The complementary subword $\bar{\beta}$ is obtained from $\alpha$ by removing all the letters that are in $\beta$. The empty word and $\alpha$ itself are also subwords of $\alpha$. For instance if $\alpha = [3,2,2,1,4,5,2,1,4]$ than $\beta = [2,1,4]$ is a subword of $\alpha$ (in several different ways). Taking its manifestation as the word formed by the second, fourth and last letter of $\alpha$ the complementary subword is $\bar{\beta} = [3,2,4,5,2,1]$. There are precisely $2^m$ subwords of $\alpha$ where $m$ is the length of $\alpha$. Using this notation the comultiplication of *LieHopf* can be more explicitly described as

$$\mu(U_\alpha) = \sum_{\beta} U_\beta \otimes U_{\bar{\beta}} \tag{3.4}$$

where the sum is over all subwords $\beta$ of $\alpha$ (counting multiplicities). For example

$$\mu(U_{[1,1,3]}) = 1\otimes U_{[1,1,3]} + 2U_{[1]}\otimes U_{[1,3]} + U_{[3]}\otimes U_{[1,1]}$$
$$+U_{[1,1]}\otimes U_{[3]} + 2U_{[1,3]}\otimes U_{[1]} + U_{[1,1,3]}\otimes 1$$

Formula (3.2) says that the indeterminates $U_n$ are all primitives of *LieHopf*. And in fact the Lie algebra of primitives Prim(*LieHopf*) of the Lie Hopf algebra is the free Lie algebra on the symbols $U_1, U_2, \cdots$. Inversely the universal envelopping algebra of this Lie algebra is *LieHopf*.[1] These facts make *LieHopf* a very important object in mathematics and it has been (very) deeply studied, see especially [12].

The graded dual of the Hopf algebra *LieHopf* is a graded Hopf algebra called *Shuffle*. Here is a brief description. As an Abelian group is the free graded Abelian group with as basis all words over the *natural* numbers, i.e. the elements of the monoid of words $\mathbf{N}^*$. The duality pairing between *Shuffle* and *LieHopf* is given by

$$LieHopf \times Shuffle \longrightarrow \mathbf{Z},\ \ \langle U_\alpha, \beta\rangle = \delta_\alpha^\beta \quad \text{(Kronecker delta)} \tag{3.5}$$

The multiplication and comultiplication on *Shuffle* are determined by the requirements

$$\langle \alpha\otimes\alpha', \mu_{Sh}\rangle = \langle m_{LH}(\alpha\otimes\alpha'), \beta\rangle,\ \ \langle \alpha, m_{Sh}(\beta\otimes\beta')\rangle = \langle \mu_{LH}(\alpha), \beta\otimes\beta'\rangle \tag{3.6}$$

where the pairing on the tensor products is the obvious one

$$(LieHopf\otimes LieHopf)\times(Shuffle\otimes Shuffle)\longrightarrow \mathbf{Z},\ \ \langle U_\alpha\otimes U_{\alpha'}, \beta\otimes\beta'\rangle = \delta_\alpha^\beta\delta_{\alpha'}^{\beta'}$$

It is not difficult to see that this works out as follows: the multiplication is the shuffle product (see section 1 above) and the comultiplication is ‘cut’. Thus, if $\beta = [b_1, \cdots, b_n]$

$$m_{Sh}(\beta\otimes\beta') = \beta\times_{sh}\beta',\ \ \mu_{Sh}(\beta) = \sum_{i=0}^{n}[b_1,\cdots,b_i]\otimes[b_{i+1},\cdots,b_n]$$

The unit element is the empty word, and the counit is given by $\varepsilon([\ ]) = 1$ and $\varepsilon(\beta) = 0$ when $\beta$ is of length greater or equal to one. The grading is given by $\deg(\beta) = \mathrm{wt}(\beta)$. As always in a connected graded situation, there exists a unique antipode, see e.g. [9] or [14].

## 4. A (seemingly irrelevant) digression: endomorphisms of finite rank Hopf algebras.

In this section $H$ is a finite rank Hopf algebra (over $\mathbf{Z}$), meaning that the underlying Abelian group is free of finite rank. Let $\mathrm{End}(H)$ denote the Abelian group of Abelian group endomorphisms of $H$.

For any Abelian group $H$ there are natural morphisms

$$\begin{aligned} &H\otimes H^* \longrightarrow \mathrm{End}(H),\ \ u\otimes f\mapsto\varphi,\ \ \varphi(v) = f(v)u \\ &\mathrm{End}(H)\otimes\mathrm{End}(H)\longrightarrow\mathrm{End}(H\otimes H),\ (f\otimes g)(u\otimes v) = f(u)\otimes g(v)\end{aligned} \tag{4.1}$$

If $H$ is a free Abelian group of finite rank these two Abelian group endomorphisms are isomorphisms. (But if $H$ is of infinite rank this is definitely not the case.)

Now, as the tensor product of two Hopf algebras $H\otimes H^*$ carries a Hopf algebra structure. Here

[1] This is why I like to use *LieHopf* to designate this Hopf algebra.

is how that structure looks in terms of endomorphisms, i.e. on $\mathrm{End}(H)$.

*Convolution.* The multiplication is convolution. Let $f$ and $g$ be two elements of $\mathrm{End}(H)$. Then their product is the composite

$$H \xrightarrow{\mu_H} H \otimes H \xrightarrow{f \otimes g} H \otimes H \xrightarrow{m_H} H \tag{4.2}$$

*Coconvolution.* The comultiplication is given as follows. Let $f$ be an element of $\mathrm{End}(H)$. Then $\mu(f)$ is the following element of $\mathrm{End}(H) \otimes \mathrm{End}(H)$

$$H \otimes H \xrightarrow{m_H} H \xrightarrow{f} H \xrightarrow{\mu_H} H \otimes H \tag{4.3}$$

Of course (4.3) defines an element of $\mathrm{End}(H \otimes H)$; but then using the second of the isomorphisms of (4.1) this yields an element of $\mathrm{End}(H) \otimes \mathrm{End}(H)$.

*Unit.* The unit of $\mathrm{End}(H)$ is the endomorphism $H \xrightarrow{\varepsilon} \mathbf{Z} \xrightarrow{e} H$ of $H$.

*Counit.* The counit of $\mathrm{End}(H)$ takes an endomorphism $f$ in $\mathrm{End}(H)$ to the number $(\varepsilon \circ f \circ e)(1)$.

*Antipode.* Finally the antipode of $\mathrm{End}(H)$ is given by $f \mapsto \iota \circ f \circ \iota$.

The Hopf algebra $\mathrm{End}(H)$ is also selfdual. This is easiest seen in its guise $H \otimes H^*$. Indeed

$$(H \otimes H^*)^* \xrightarrow{\sim} H^* \otimes H^{**} \xleftarrow{\sim} H^* \otimes H \xrightarrow{\sim} H \otimes H^*$$

The corresponding inner product is

$$(H \otimes H^*) \times (H \otimes H^*) \longrightarrow \mathbf{Z}, \ \ \langle u \otimes f, v \otimes g \rangle = f(v) g(u)$$

which is nondegenerate but not positive definite.

At the level of $\mathrm{End}(H)$ the autoduality is a combination of a canonical pairing

$$\mathrm{End}(H) \otimes \mathrm{End}(H^*) \longrightarrow \mathbf{Z} \tag{4.4}$$

and the isomorphism $\mathrm{End}(H) \longrightarrow \mathrm{End}(H^*)$ that assigns to an endomorphism of $H$ the dual endomorphism which is an endomorphism of $H^*$. The canonical pairing (4.4) is defined as follows. Take a basis $u_1, \cdots, u_n$ of $H$ and let $v^1, \cdots, v^n$ be the dual basis of $H^*$ (so that $v^i(u_j) = \delta^i_j$ (Kronecker delta). Then

$$\gamma = \sum_{i=1}^{n} u_i \otimes v^i \ \in \ H \otimes H^* \tag{4.5}$$

is a special element of $H \otimes H^*$ that is independent of the choice of basis. The easiest way to see this is to remark that the element $\gamma$ is the image of $1 \in \mathbf{Z}$ under the dual of the evaluation morphism

$$\mathrm{ev}: \ H \otimes H^* \longrightarrow \mathbf{Z}, \ \ u \otimes f \mapsto f(u) \tag{4.6}$$

The pairing (4.4) is now defined by

$$f \otimes g \mapsto (f \otimes g)(\gamma) = \sum_{i=1}^{n} f(u_i)g(v^i) \tag{4.7}$$

The elements of $\mathrm{End}(H)$ are endomorphisms. So there is a second multiplication on $\mathrm{End}(H)$, viz composition. This second multiplication is not necessarily distributive over the first (which would make $\mathrm{End}(H)$ a ring object in the category of coalgebras. Still an extra multiplication, i.e a second way of producing a new element from two given ones, can be very useful even when it has no particular compatibility properties with respect to the other structure present. Examples of this are the many extra 'multiplications' of divided power sequences as they are used in [6] to describe a basis over the integers of the Lie algebra of primitives $\mathrm{Prim}(NSymm)$.

Dually there is also a second comultiplication (cocomposition).

4.8. *Open problem.* Which Hopf algebras, necessarily of square rank, are of the form $\mathrm{End}(H)$?

More generally there is a Hopf algebra structure on $\mathrm{Mod}(H,K)$ where $H$ and $K$ are two possibly different Hopf algebras with free finite rank underlying Abelian groups and there is the same open problem vis à vis these Hopf algebras. The same open problems can be considered over other base rings, for instance fields.

Note that practically everything above works perfectly fine for graded Hopf algebras (with graded dual replacing dual). The sole exception is the comultiplication (given by coconvolution). But that is is formidable exception and obstruction.

The seeming irrelevancy of these considerations lies in the following. Of course for a graded Hopf algebra $H \otimes H^{\mathrm{gr}*}$ is again a Hopf algebra; but it is a very small part of the Abelian group of homogenous endomorphisms of $H$. Indeed it is easy to check that a homogenous endomorphism of $H$ is in (the image of) $H \otimes H^{\mathrm{gr}*}$ (in $\mathrm{End}(H)$) iff it has a finite rank image.[2]

Now, as remarked before, the elements of $MPR$ (i.e. $\mathbf{Z}$-linear combinations of permutations) are to be interpreted as homogeneous endomorphisms of the Hopf algebra *Shuffle*. However, as it will turn out, with the exception of the scalar multiples of the empty permutation, none of these endomorphisms lies in $H \otimes H^{\mathrm{gr}*}$.

## 5. *MPR* as a Hopf algebra of endomorphisms.

The first thing to do is to interprete ekements of $MPR$ as endomorphisms (of Abelian groups) of *Shuffle.*

*Permutations as Shuffle endomorphisms.* Let $\sigma = [s_1, \cdots, s_m]$ be a permutation word. The corresponding permutation is of course $j \mapsto s_j$. Further let $\alpha$ be a basis element, i.e. a word, of *Shuffle.* Then $\sigma$ takes $\alpha$ to zero unless $\lg(\alpha) = m = \lg(\sigma)$ and if $\alpha = [a_1, \cdots, a_m]$ is of length $m$ then

$$\sigma(\alpha) = [a_{\sigma(1)}, \cdots, a_{\sigma(m)}] = [a_{s_1}, \cdots, a_{s_m}] \tag{5.1}$$

Note that unless $\sigma$ is the empty permutation (which takes [ ] to [ ] and is zero on all other words) $\sigma$ is never in $Shuffle \otimes Shuffle^{gr*}$ (because for every length $\geq 1$ there are infinitely many words of that length).

[2] More generally for an infinite rank free Abelian group $M$ or vectorspace an endomorphism is in $M \otimes M^*$ iff the images of both the morphism itself and its dual are of finite rank.

*Convolution.* Here is how convolution of two permutations works out with this interpretation of permutations as endomorphisms of *Shuffle*. So let $\sigma = [s_1, \cdots, s_m]$ and $\tau = [t_1, \cdots, t_n]$ be two permutation words. Their convolution[3] is given by

$$Shuffle \xrightarrow{\mu_{Sh}} Shuffle \otimes Shuffle \xrightarrow{\sigma \otimes \tau} Shuffle \otimes Shuffle \xrightarrow{m_{Sh}} Shuffle$$

Because $\mu_{Sh}$ is cut, the middle morphism is zero on all terms of $\mu_{Sh}(\alpha)$ for all words $\alpha$ which are not of length $m+n$. And if $\alpha = [a_1, \cdots, a_m, a_{m+1}, \cdots, a_{m+n}]$, $\sigma \otimes \tau$ is zero on all terms of $\mu_{Sh}(\alpha)$ except the summand $[a_1, \cdots, a_m] \otimes [a_{m+1}, \cdots, a_{m+n}]$ and this summand is taken to $[a_{s_1}, \cdots, a_{s_m}] \otimes [a_{m+t_1}, \cdots, a_{m+t_n}]$. And thus the convolution product of $\sigma$ and $\tau$ takes $\alpha$ to the shuffle product $[a_{s_1}, \cdots, a_{s_m}] \times_{sh} [a_{m+t_1}, \cdots, a_{m+t_n}]$ and it follows immediately that the convolution of $\sigma$ and $\tau$ is equal to to the sum of permutation words $[s_1, \cdots, s_m] \times_{sh} [m+t_1, \cdots m+t_n]$ which is the multiplication on $MPR$ (see section 1 above). There is nothing new about this; this is (more or less[4]) the way the multiplication of $MPR$ was introduced in [10].

*Coconvolution.* The coconvolution construction takes an endomorphism $f$ of *Shuffle* to the composite morphism

$$Shuffle \otimes Shuffle \xrightarrow{m_{Sh}} Shuffle \xrightarrow{f} Shuffle \xrightarrow{\mu_{Sh}} Shuffle \otimes Shuffle$$

that is, it defines a morphism of Abelian groups

$$MPR \subset \mathrm{End}(Shuffle) \xrightarrow{coconv} \mathrm{End}(Shuffle \otimes Shuffle) \tag{5.2}$$

And to turn this into a comultiplication on $MPR$ some sort of projection is needed from the image of coconv in $\mathrm{End}(Shuffle \otimes Shuffle)$ to $MPR \otimes MPR \subset \mathrm{End}(Shuffle) \otimes \mathrm{End}(Shuffle)$. Here is an example how that might be accomplished.

Take $\sigma = [3,1,4,5,2]$ and lets see what the coconvolution of this does to an element of the form $[a_1, a_2] \otimes [b_3, b_4, b_5]$ from $Shuffle \otimes Shuffle$. (Different letters are used for the two copies of *Shuffle* to make it easier to identify where each is from.)

$$\begin{aligned}
&[a_1,a_2]\otimes[b_3,b_4,b_5] \underset{\times_{sh}}{\longmapsto} [a_1,a_2,b_3,b_4,b_5]+[a_1,b_3,a_2,b_4,b_5]+\\
&+[a_1,b_3,b_4,a_2,b_5]+[a_1,b_3,b_4,b_5,a_2]+[b_3,a_1,a_2,b_4,b_5]+[b_3,a_1,b_4,a_2,b_5]+\\
&+[b_3,a_1,b_4,b_5,a_2]+[b_3,b_4,a_1,a_2,b_5]+[b_3,b_4,a_1,b_5,a_2]+[b_3,b_4,b_5,a_1,a_2]\\
&\underset{\sigma}{\longmapsto} [b_3,a_1,b_4,b_5,a_2]+[a_2,a_1,b_4,b_5,b_3]+\cdots+[b_5,b_3,a_1,a_2,b_4]
\end{aligned}$$

and now cut has to be applied to the final ten terms. But we are looking for an endomorphisms of the form $\in MPR \otimes MPR$. So the only cuts that can contribute are ones which have only $a$'s on the left and only $b$'s on the right. The only one of the final ten terms for which this is possible is the second one and there is (of course) only one cut of this term which qualifies, yielding

[3] Quite generally if $C$ is a coalgebra and $A$ is an algebra, the convolution of two morphisms $f, g\colon C \longrightarrow A$ is the morphism $C \xrightarrow{\mu_C} C \otimes C \xrightarrow{f \otimes g} A \otimes A \xrightarrow{m_A} A$. For the right Hopf algebras this is indeed the classical convolution of functions, see [9] or [S. Dascalescu, C. Natascescu and S. Raianu, *Hopf algebras. An introduction*, Marcel Dekker, 2001].

[4] Actually in [10] permutations are seen as endomorphisms of the underlying group of *LieHopf* and it is the dual multiplication on $MPR$, see section 9, that is defined by convolution.

$$[a_2,a_1]\otimes[b_4,b_5,b_3] = (\mathrm{st}([3,1])\otimes\mathrm{st}([4,5,2]))([a_1,a_2]\otimes[b_3,b_4,b_5])$$

And thus according to this procedure the $(\lg = 2)\otimes(\lg = 3)$ component of $\mu_{MPR}([3,1,4,5,2])$ is equal to $\mathrm{st}([3,1])\otimes\mathrm{st}([4,5,2])$ exactly as it should be according to the description given in section 1 above.

This works in general and gives the correct description of the comultiplication on $MPR$.

All the same this is a very dodgy procedure. For instance there is a very nice generalization of the way $MPR$ acts as endomorphisms on *Shuffle* for arbitrary words. This goes as follows. Let now $\sigma = [s_1,\cdots s_n]$ be an arbitrary word of height $m = ht(\sigma) = \max\{s_1,\cdots,s_n\}$. Then $\sigma$ acts as zero on all basis words $\alpha$ of *Shuffle* of length $\neq m$ and if $\alpha = [a_1,\cdots a_m]$ is a word of length $m$ $\sigma(\alpha) = [a_{s_1},\cdots,a_{s_n}]$.

Now apply the procedure used above for $MPR$. This does define a comultiplication, but the Hopf property (viz. multiplication is a coalgebra morphism, or, equivalently, the comultiplication is an algebra morphism [5]) fails completely.[6]

However, as it turns out there is a rather different way to make arbitrary words act on *Shuffle* which does yield a Hopf algebra and for which the dodgy procedure outlined above in the case of $MPR$, does work.

What is needed to use coconvolution to define a Hopf algebra on a suitable module of endomorphisms is some kind of suitable projection $\mathrm{End}(H\otimes H)\xrightarrow{\pi}\mathrm{End}(H)\otimes\mathrm{End}(H)$ (where $\mathrm{End}(H)\otimes\mathrm{End}(H)$ is seen as a submodule of $\mathrm{End}(H\otimes H)$ in the natural way). The next section is devoted to a preliminary analysis of the question of what projections could work.

## 6. Hopf algebras of endomorphisms. The general problem.

So, now let $E(H)\subset\mathrm{End}(H)$ be a submodule that is closed under convolution; i.e. if $f,g\in E(H)$ then so is $m\circ(f\otimes g)\circ\mu$.

As remarked before, coconvolution, $f\mapsto\mu\circ f\circ m$, takes an endomorphism $f$ to an element in $\mathrm{End}(H\otimes H)$ and not necessarily (or even usually) to an element of $E(H)\otimes E(H)$.

In this connection note that if $H$ is of infinite rank (dimension) then $\mathrm{End}(H)\otimes\mathrm{End}(H)$ is only a small part of $\mathrm{End}(H\otimes H)$. Indeed, let $\{e_i\}_{i\in I}$ be a basis of $H$. Then $\{e_i\otimes e_j\}_{i,j\in I}$ is a basis of $H\otimes H$ and an endomorphism $f$ of $H\otimes H$ is given by an array of coefficients $(c_{i,j}^{r,s})$

$$f(e_i\otimes e_j) = \sum^{<\infty} c_{i,j}^{r,s} e_r\otimes e_s$$

where for each $(i,j)$ there are only finitely many $(r,s)$ for which $c_{i,j}^{r,s}$ is nonzero. But there is no bound on how many of these coefficients are nonzero and that means that the matrix $(c_{i,j}^{r,s})$ with columns indexed by pairs $\begin{pmatrix} s \\ j \end{pmatrix}$ and rows indexed by pairs $\begin{pmatrix} r \\ i \end{pmatrix}$ is usually of infinite rank.

Now let $f_k, g_k$ be endomorphism of $H$ given by coefficients $(a_{i,k}^r)$ and $(b_{j,k}^s)$, $k = 1,\cdots,t$, so that

[5] This is most always the hardest property to check.

[6] Some two years ago I wasted a couple of months research time to try to fix things up; typical Ptolemaic-epicycle-type thinking, and thoroughly useless.

$$f_k(e_i) = \sum a_{i,k}^r e_r, \quad g_k(e_j) = \sum b_{j,k}^s e_s$$

It follows that for $f = \sum_{k=1}^{t} g_k \otimes h_k$ to hold one must have the matrix identity

$$\begin{pmatrix} \vdots & \vdots & \cdots & \vdots \\ a_{i,1}^r & a_{i,2}^r & \cdots & a_{i,t}^r \\ \vdots & \vdots & \cdots & \vdots \end{pmatrix} \begin{pmatrix} \cdots & b_{j,1}^s & \cdots \\ \cdots & b_{j,2}^s & \cdots \\ \vdots & \vdots & \vdots \\ \cdots & b_{j,t}^s & \cdots \end{pmatrix} = \begin{pmatrix} \cdots & \vdots & \cdots \\ \cdots & c_{i,j}^{r,s} & \cdots \\ \cdots & \vdots & \cdots \end{pmatrix}$$

which is usually impossible to solve because, generically, the matrix on the right hand side has infinite rank, while the product on the left has at most rank $t$.

Things become even worse if in the graded case one only looks at homogeneous endomorphisms. In that case a sum $f = \sum_{k=1}^{t} g_k \otimes h_k$ restricted to $(H \otimes H)_n$ is necessarily of block diagonal form corresponding to the decomposition

$$(H \otimes H)_n = \bigoplus_{i=0}^{n} H_i \otimes H_{n-i}$$

which of course need not be the case for an arbitrary homogeneous endomorphism of $H \otimes H$.[7]

So what is needed for coconvolution $\mathrm{End}(H) \longrightarrow \mathrm{End}(H \otimes H)$ to induce a comultiplication is some kind of nice projection $\mathrm{End}(H \otimes H) \xrightarrow{\pi} E(H) \otimes E(H)$ so that the comultiplication on $E(H)$ is defined by

$$E(H) \xrightarrow{\mu_{E(H)}} E(H) \otimes E(H), \quad f \mapsto \pi \circ \mu_H \circ f \circ m_H \tag{6.1}$$

This comultiplication needs to be coassociative which is not automatic and poses some condition on $\pi$. However, as always, the main problem is to guarantee the Hopf property (the comultiplication is an algebra morphism, or equivalently, the multiplication is a coalgebra morphism). This means that the following diagram must be commutative.

$$\begin{array}{ccccc} E(H)^{\otimes 2} & \xrightarrow{coconv^{\otimes 2}} & \mathrm{End}(H^{\otimes 2}) \otimes \mathrm{End}(H^{\otimes 2}) & \xrightarrow{\pi \otimes \pi} & E(H)^{\otimes 4} \\ \downarrow \mathrm{conv} & & \downarrow \subset & & \downarrow \mathrm{id} \otimes \mathrm{tw} \otimes \mathrm{id} \\ \mathrm{E(H)} & & \mathrm{End}(H^{\otimes 4}) & & E(H)^{\otimes 4} \\ \downarrow \mathrm{coconv} & & \downarrow \alpha & & \downarrow \mathrm{conv}^{\otimes 2} \\ \mathrm{End}(H^{\otimes 2}) & & \mathrm{End}(H^{\otimes 4}) & & E(H)^{\otimes 2} \\ \downarrow = & & \downarrow \mathrm{conv}^{(2)} & & \downarrow = \\ \mathrm{End}(H^{\otimes 2}) & \xrightarrow{=} & \mathrm{End}(H^{\otimes 2}) & \xrightarrow{\pi} & E(H)^{\otimes 2} \end{array} \tag{6.2}$$

Here $\alpha$ is conjugation by $(\mathrm{id} \otimes \mathrm{tw} \otimes \mathrm{id})$, $\mathrm{tw}(a \otimes b) = b \otimes a$, i.e.

$$\alpha(f) = (\mathrm{id} \otimes \mathrm{tw} \otimes \mathrm{id}) \circ f \circ (\mathrm{id} \otimes \mathrm{tw} \otimes \mathrm{id})$$

[7] It is not even true that all block diagonal endomorphims are in $\mathrm{End}(H) \otimes \mathrm{End}(H)$ even when it is assumed that there is a bound on the rank of the homogeneous components.

That this is the right ‘twist’ to use follows from the commutativity of the diagram

$$\begin{array}{ccc} \mathrm{End}(H^{\otimes 4}) & \xleftarrow{\supset} & \mathrm{End}(H)^{\otimes 4} \\ \downarrow \alpha & & \downarrow \mathrm{id}\otimes \mathrm{tw}_{\mathrm{End}(H)}\otimes \mathrm{id} \\ \mathrm{End}(H^{\otimes 4}) & \xleftarrow{\supset} & \mathrm{End}(H)^{\otimes 4} \end{array}$$

Indeed,

$$\begin{aligned} &(\mathrm{id}\otimes \mathrm{tw}\otimes \mathrm{id})(f_1\otimes f_2\otimes f_3\otimes f_4)(\mathrm{id}\otimes \mathrm{tw}\otimes \mathrm{id})(a_1\otimes a_2\otimes a_3\otimes a_4) = \\ &(\mathrm{id}\otimes \mathrm{tw}\otimes \mathrm{id})(f_1(a_1)\otimes f_2(a_3)\otimes f_3(a_2)\otimes f_4(a_4)) = \\ &= f_1(a_1)\otimes f_3(a_2)\otimes f_2(a_3)\otimes f_4(a_4) \end{aligned}$$

while

$$\begin{aligned} &((\mathrm{id}\otimes \mathrm{tw}_{\mathrm{End}(H)}\otimes \mathrm{id})(f_1\otimes f_2\otimes f_3\otimes f_4))(a_a\otimes a_2\otimes a_3\otimes a_4) = \\ &(f_1\otimes f_3\otimes f_2\otimes f_4)(a_a\otimes a_2\otimes a_3\otimes a_4) = \\ &= f_1(a_1)\otimes f_3(a_2)\otimes f_2(a_3)\otimes f_4(a_4) \end{aligned}$$

The morphism marked ‘$\mathrm{conv}^{(2)}$‘ in diagram (6.2) above is given by

$$conv^{(2)}(g) = (m_H\otimes m_H)g(\mu_H\otimes \mu_H)$$

The left half of diagram (6.2) is always commutative. This is part of a Hopf-like structure on the collection of modules

$$\mathrm{End}(H^{\otimes n}),\ n = 1,2,\cdots$$

defined by convolution and coconvolution. Here convolution is seen as the morphism

$$\mathrm{End}(H^{\otimes 2}) \longrightarrow \mathrm{End}(H),\quad f \mapsto m_H \circ f \circ \mu_H \tag{6.3}$$

and, as before, coconvolution is the morphism

$$\mathrm{End}(H) \longrightarrow \mathrm{End}(H^{\otimes 2}),\quad f \mapsto \mu_H \circ f \circ m_H \tag{6.4}$$

More generally there are $n$ morphisms

$$\begin{aligned} &\mathrm{End}(H^{\otimes n+1}) \xrightarrow{\mathrm{conv}_{i,i+1}} \mathrm{End}(H^n),\quad i = 1,\cdots,n \\ &f \mapsto (\underbrace{\mathrm{id}\otimes\cdots\otimes \mathrm{id}}_{i-1}\otimes m_H\otimes \underbrace{\mathrm{id}\otimes\cdots\otimes\mathrm{id}}_{n-i})\circ f\circ(\underbrace{\mathrm{id}\otimes\cdots\otimes \mathrm{id}}_{i-1}\otimes \mu_H\otimes \underbrace{\mathrm{id}\otimes\cdots\otimes\mathrm{id}}_{n-i}) \end{aligned}$$

and $n$ morphisms

$$\begin{aligned} &\mathrm{End}(H^{\otimes n}) \xrightarrow{\mathrm{coconv}_{i\,i+1}} \mathrm{End}(H^{n+1}),\quad i = 1,\cdots,n \\ &f \mapsto (\underbrace{\mathrm{id}\otimes\cdots\otimes \mathrm{id}}_{i-1}\otimes \mu_H\otimes \underbrace{\mathrm{id}\otimes\cdots\otimes\mathrm{id}}_{n-i})\circ f\circ(\underbrace{\mathrm{id}\otimes\cdots\otimes \mathrm{id}}_{i-1}\otimes m_H\otimes \underbrace{\mathrm{id}\otimes\cdots\otimes\mathrm{id}}_{n-i}) \end{aligned}$$

and these ‘mutiplication-like’ and ‘comultiplication-like’ structures are associative and coassociatative in the sense that

$$\mathrm{conv}\circ\mathrm{conv}_{1,2}=\mathrm{conv}\circ\mathrm{conv}_{2,3}:\ \mathrm{End}(H^{\otimes 3})\longrightarrow\mathrm{End}(H)$$

$$\mathrm{coconv}_{1,2}\circ\mathrm{coconv}=\mathrm{coconv}_{2,3}\circ\mathrm{coconv}:\ \mathrm{End}(H)\longrightarrow\mathrm{End}(H^{\otimes 3})$$

and there is a Hopf-like property in that the following diagram is commutative

$$\begin{array}{ccc}
\mathrm{End}(H^{\otimes 2}) & \xrightarrow{\mathrm{coconv}_{H^{\otimes 2}}} & \mathrm{End}(H^{\otimes 4}) \\
\downarrow \mathrm{conv}_H & & \downarrow \alpha \\
\mathrm{End}(H) & & \mathrm{End}(H^{\otimes 4}) \\
\downarrow \mathrm{coconv} & & \downarrow \mathrm{conv}^{(2)} \\
\mathrm{End}(H^{\otimes 2}) & = & \mathrm{End}(H^{\otimes 2})
\end{array} \qquad (6.5)$$

In turn, the commutativity of diagram (6.5) follows from the Hopf property of the Hopf algebra $H$ itself in that the following diagram is commutative

$$\begin{array}{ccccc}
H^{\otimes 2} & = & H^{\otimes 2} & & \\
\downarrow \mu^{\otimes 2} & & \downarrow m & & \\
H^{\otimes 4} & & H & & \\
\downarrow \mathrm{id}\otimes\mathrm{tw}\otimes\mathrm{id} & & \downarrow \mu & & \\
H^{\otimes 4} & \xrightarrow{m^{\otimes 2}} & H^{\otimes 2} & & \\
 & & \downarrow f & & \\
 & & H^{\otimes 2} & \xrightarrow{\mu^{\otimes 2}} & H^{\otimes 4} \\
 & & \downarrow m & & \downarrow \mathrm{id}\otimes\mathrm{tw}\otimes\mathrm{id} \\
 & & H & & H^{\otimes 4} \\
 & & \downarrow m & & \downarrow m^{\otimes 2} \\
 & & H^{\otimes 2} & = & H^{\otimes 2}
\end{array}$$

This Hopf-like structure on the collection of modules $\mathrm{End}(H^{\otimes n}),\ n=1,2,\cdots$ deserves, I think, a great deal more investigation.

6.6. *Open problem*. What kind of properties should the projection $\pi$ have in order that the right hand side of diagram (6.2) be commutative.

More generally one can consider submodules $E(H)$ that are not necessarily stable under convolution (together with a second projection $\mathrm{End}(H)\longrightarrow E(H)$ and instead of with endomorphisms one can work with module morphisms $\mathrm{End}(H^{\otimes n},K^{\otimes n})$ where $K$ is a second Hopf algebra.

**7. The double word Hopf algebra** *dWHA*.

Let $\mathcal{X}=\{x_1,x_2,\cdots\}$ be an auxiliary alphabet. A basis as an Abelian group of *dWHA* is formed by pairs of words in the auxiliary alphabet $\mathcal{X}$ with equal support

$$p=\begin{pmatrix}\rho\\ \sigma\end{pmatrix},\ \ \mathrm{supp}(\rho)=\mathrm{supp}(\sigma)$$

Here the actual symbols that occur are not important; it is only the pattern of $\rho$ and $\sigma$ relative to each other that is relevant. Thus for instance

$$\begin{pmatrix}[x_1,x_2,x_1,x_3,x_3,x_1,x_4]\\ [x_2,x_3,x_2,x_4,x_1]\end{pmatrix}, \begin{pmatrix}[x_7,x_6,x_7,x_2,x_2,x_7,x_5]\\ [x_6,x_2,x_6,x_5,x_7]\end{pmatrix}, \begin{pmatrix}[y_3,z_4,y_3,x_2,x_2,y_3,x_1]\\ [z_4,x_2,z_4,x_1,y_3]\end{pmatrix} \tag{7.1}$$

all denote the same basis element of *dWHA*. I call these things substitutions.

These elements $p$ can (and often should) be thought of as defining endomorphisms of *Shuffle*; more precisely they are recipes for endomorphisms as follows. The 'substitution' $p$ acts as zero on all words over **N**, i.e. the canonical basis elements of *Shuffle,* that are not of the same pattern as $\rho$ and if it is of the same pattern as $\rho$ then it is taken into the corresponding basis element of *Shuffle* represented by the pattern $\sigma$.

Thus, for example, if $p$ is the substitution (7.1) and $\alpha = [a_1, a_2, \cdots, a_m]$

$$p(\alpha) = \begin{cases} 0 \text{ unless } \lg(\alpha) = m = 7 \text{ and } a_1 = a_3 = a_6,\ a_4 = a_5 \\ [a_2, a_3, a_2, a_4, a_1] \text{ if } \lg(\alpha) = m = 7 \text{ and } a_1 = a_3 = a_6,\ a_4 = a_5 \end{cases}$$

Obviously these endomorphisms satisfy a 'homogeneity' property; they act the same everywhere. For instance if $\phi : \mathbf{N} \longrightarrow \mathbf{N}$ is any map and $\phi_* : \mathbf{N}^* \longrightarrow \mathbf{N}^*$ denotes the corresponding induced map on words $\phi_*(\alpha) = [\phi(a_1), \phi(a_2), \cdots, \phi(a_m)]$

$$\phi_* \circ p = p \circ \phi_* \tag{7.2}$$

7.3. *Open problem.* Characterize the recipe endomorphisms $p$ more precisely.

The next step is to describe the graded Hopf algebra structure on *dWHA*.

*Underlying Abelian group.* The uderlying Abelian group is the countable free Abelian group with as basis all substitutions $p$. Note that this not a basis of the tensor product of Hopf algebras *LieHopf* $\otimes$ *Shuffle* but a suitable quotient set obtained by identifying substitutions of the same patterns. Included is the 'empty substitution'

$$\begin{pmatrix}[\ ]\\ [\ ]\end{pmatrix}$$

which acts on *Shuffle* by taking the empty word [ ] to itsef and every other basis element of *Shuffle* to zero.

*Grading.* The grading on *dWHA* is given by

$$\deg(p) = \#\mathrm{supp}(\rho) \tag{7.4}$$

For example the degree of the basis element (7.1) is 4. The degree of the empty substitution is zero and that is the only basis element of degree zero so that the graded Abelian group *dWHA* is connected. Note that the rank of each homogeneous piece of *dWHA* is infinite.

*Multiplication.* Let

$$p = \begin{pmatrix}\rho\\ \sigma\end{pmatrix},\quad p' = \begin{pmatrix}\rho'\\ \sigma'\end{pmatrix} \tag{7.5}$$

be two substitutions. If necessary, first rewrite the second one (or the first one, or both) so that $\text{supp}(\rho) \cap \text{supp}(\rho') = \emptyset$. Then the product of the two substitutions (7.4) is the sum of substitutions

$$m_{dWHA}(p \otimes p') = \begin{pmatrix} \rho * \rho' \\ \sigma \times_{sh} \sigma \end{pmatrix} \tag{7.6}$$

where $*$ denotes concatenation, $\times_{sh}$ is the shuffle product, and if $u = \sigma_1 + \sigma_2 + \cdots + \sigma_r$ is a sum of words with $\text{supp}(\sigma_1) = \text{supp}(\sigma_2) = \cdots = \text{supp}(\sigma_r) = \text{supp}(\rho)$

$$\begin{pmatrix} \rho \\ u \end{pmatrix} = \begin{pmatrix} \rho \\ \sigma_1 \end{pmatrix} + \begin{pmatrix} \rho \\ \sigma_2 \end{pmatrix} + \cdots + \begin{pmatrix} \rho \\ \sigma_r \end{pmatrix}$$

*Unit element.* The unit element is the empty substitution.

It is easy to see that the multiplication is associative (and that the empty substitution indeed acts as a unit element. Further, clearly the multiplication respects the grading making $(dWHA, m, e)$ a connected graded algebra.

*Comultiplication.* To write down the comultiplication of $dWHA$ a preliminary definition is needed. Let

$$\alpha = [a_1, \cdots, a_m]$$

be a word over an alphabeth $\mathcal{X}$. A good cut of $\alpha$ is a cut $[a_1, \cdots, a_r] \otimes [a_{r+1}, \cdots a_m]$ such that $\text{supp}([a_1, \cdots, a_r]) \cap \text{supp}([a_{r+1}, \cdots, a_m]) = \emptyset$. The two trivial cuts are always good. For example the good cuts of $[x_2, x_3, x_2, x_4, x_1]$ are $1 \otimes [x_2, x_3, x_2, x_4, x_1]$, $[x_2, x_3, x_2] \otimes [x_4, x_1]$, $[x_2, x_3, x_2, x_4] \otimes [x_1]$ and $[x_2, x_3, x_2, x_4, x_1] \otimes 1$ (where as usual $1$ is short for $[\,]$). A subword of $\alpha$ is a word $[a_{i_1}, a_{i_2}, \cdots, a_{i_r}]$ with $i_1 < i_2 < \cdots < i_r$.

The comultiplication of $dWHA$ is now

$$\mu_{dWHA}(p) = \sum_{\sigma_1 * \sigma_2 = \sigma} \begin{pmatrix} p^{-1}(\sigma_1) \\ \sigma_1 \end{pmatrix} \otimes \begin{pmatrix} p^{-1}(\sigma_2) \\ \sigma_2 \end{pmatrix} \tag{7.7}$$

where the sum is over all good cuts $\sigma_1 * \sigma_2 = \sigma$ of the word $\sigma$. Here $p^{-1}(\sigma_i)$ is the unique maximal subword of $\rho$ with the same support as $\sigma_i$. This respects the grading and the support condition. For instance the comultiplication of the substitution

$$p = \begin{pmatrix} [x_1, x_2, x_1, x_3, x_3, x_1, x_4, x_1, x_4] \\ [x_2, x_3, x_2, x_4, x_1] \end{pmatrix}$$

is

$$\mu_{dWHA}(p) = 1 \otimes p + \begin{pmatrix} [x_2, x_3, x_3] \\ [x_2, x_3, x_2] \end{pmatrix} \otimes \begin{pmatrix} [x_1, x_1, x_1, x_4, x_1, x_4] \\ [x_4, x_1] \end{pmatrix}$$
$$+ \begin{pmatrix} [x_2, x_3, x_3, x_4, x_4] \\ [x_2, x_3, x_2, x_4] \end{pmatrix} \otimes \begin{pmatrix} [x_1, x_1, x_1, x_1] \\ [x_1] \end{pmatrix} + p \otimes 1$$

*Counit.* The counit is given by $\varepsilon(p)=0$ if $\deg(p)>0$ and $\varepsilon$ takes the value 1 on the empty substitution.

It is now easy to check that $(dWHA,\mu,\varepsilon)$ is a coassociative connected graded coalgebra.

7.8. *Theorem.* The Hopf property holds; i.e $(dWHA,m,\mu,e,\varepsilon)$ is a connected graded bialgebra and hence (because it is connected graded) there is also an antipode, making it a Hopf algebra.

*Proof.* To prove the Hopf property it is needed to prove the commutativity of the following diagram

$$\begin{array}{ccc} dWHA^{\otimes 2} & \xrightarrow{\mu\otimes\mu} & dWHA^{\otimes 4} \\ \downarrow m & & \downarrow \mathrm{id}\otimes\mathrm{tw}\otimes\mathrm{id} \\ dWHA & & dWHA^{\otimes 4} \\ \downarrow \mu & & \downarrow m\otimes m \\ dWHA^{\otimes 2} & = & dWHA^{\otimes 2} \end{array}$$

So let

$$p=\begin{pmatrix}\rho\\ \sigma\end{pmatrix} \quad \text{and} \quad p'=\begin{pmatrix}\rho'\\ \sigma'\end{pmatrix}$$

be two substitutions. Their product is

$$m(p\otimes p')=\begin{pmatrix}\rho*\rho'\\ \sigma\times_{sh}\sigma'\end{pmatrix}$$

Now consider a cut $\gamma=\gamma_1*\gamma_2$ of a shuffle $\gamma$ of $\sigma$ and $\sigma'$. If this is a good cut it induces good cuts of $\sigma$ and $\sigma'$, say $\sigma=\sigma_1*\sigma_2$ and $\sigma'=\sigma_1'*\sigma_2'$ because $\mathrm{supp}(\sigma)\cap\mathrm{supp}(\sigma')=\varnothing$. Indeed $\sigma_1$ is the prefix of $\sigma$ consisting of all letters of $\sigma$ that occur in the prefix $\gamma_1$ of $\gamma$. (Note that these letters are recognizable because of the support condition and that they form a prefix (not just a subword) because in a shuffle the letters of each of the two factors occur in their original order.) Also $\sigma_2$ is the suffix of $\sigma$ consisting of the letters of $\sigma$ that occur in $\gamma_2$. Moreover, $\gamma_1$ is a shuffle of $\sigma_1$ and $\sigma_1'$ and $\gamma_2$ is a shuffle of $\sigma_2$ and $\sigma_2'$. Inversely, let $\sigma=\sigma_1*\sigma_2$ and $\sigma'=\sigma_1'*\sigma_2'$ be two good cuts, let $\gamma_1$ be any shuffle of $\sigma_1,\sigma_1'$ and $\gamma_2$ a shuffle of $\sigma_2,\sigma_2'$, then $\gamma_1*\gamma_2$ is a shuffle of $\sigma,\sigma'$ and all shuffles of $\sigma,\sigma'$ are obtained this way. Thus the result of applying $\mu\circ m$ to $p\otimes p'$ is the sum

$$\sum\begin{pmatrix}p^{-1}(\sigma_1)*p'^{-1}(\sigma_1')\\ \sigma_1\times_{sh}\sigma_1'\end{pmatrix}\otimes\begin{pmatrix}p^{-1}(\sigma_2)*p'^{-1}(\sigma_2')\\ \sigma_2\times_{sh}\sigma_2'\end{pmatrix} \tag{7.9}$$

where the sum is over all good cuts $\sigma=\sigma_1*\sigma_2$, $\sigma'=\sigma_1'*\sigma_2'$. Here it is also necessary to note that if

$$q=\begin{pmatrix}\rho*\rho'\\ \gamma\end{pmatrix}, \quad \text{where } \gamma \text{ is a term from the shuffle product of } \sigma,\sigma'$$

and

$\gamma = \gamma_1 * \gamma_2$ is a good cut with $\gamma_1$ a shuffle of $\sigma_1, \sigma_1'$; $\gamma_2$ a shuffle of $\sigma_2, \sigma_2'$

then

$$q^{-1}(\gamma_1) = p^{-1}(\sigma_1) * p'^{-1}(\sigma_1'), \quad q^{-1}(\gamma_2) = p^{-1}(\sigma_2) * p'^{-1}(\sigma_2')$$

(This uses the support condition again.) However, the sum (7.9) is precisely what one gets by applying $(m \otimes m) \circ (\mathrm{id} \otimes \mathrm{tw} \otimes \mathrm{id}) \circ (\mu \otimes \mu)$ to $p \otimes p'$. This proves the theorem (modulo the trivial verifications concerning unit and counit).

There is a natural non positive definite inner product on $dWHA$ defined by

$$\left\langle \begin{pmatrix} \rho \\ \sigma \end{pmatrix}, \begin{pmatrix} \rho' \\ \sigma' \end{pmatrix} \right\rangle = \begin{cases} 1 & \text{if } \rho = \sigma',\ \rho' = \sigma \\ 0 & \text{otherwise} \end{cases} \tag{7.10}$$

where of course the two substitutions must be so written that $\mathrm{supp}(\rho) = \mathrm{supp}(\sigma) = \mathrm{supp}(\rho') = \mathrm{supp}(\sigma')$. So the more careful statement is that

$$\left\langle \begin{pmatrix} \rho \\ \sigma \end{pmatrix}, \begin{pmatrix} \rho' \\ \sigma' \end{pmatrix} \right\rangle = 1$$

iff there is a substitution of variables (of the variables in $\rho'$ and $\sigma'$) such that $\rho = \sigma',\ \rho' = \sigma$, and otherwise this inner product is zero.

7.11. *Theorem.* The Hopf algebra $dWHA$ is selfdual with respect to the inner product (7.10).

*Proof.* What needs to be shown is that

$$\left\langle \begin{pmatrix} \rho \\ \sigma \end{pmatrix} \otimes \begin{pmatrix} \rho' \\ \sigma' \end{pmatrix}, \sum \begin{pmatrix} p''^{-1}(\sigma_1'') \\ \sigma_1'' \end{pmatrix} \otimes \begin{pmatrix} p''^{-1}(\sigma_2'') \\ \sigma_2'' \end{pmatrix} \right\rangle = \left\langle \begin{pmatrix} \rho * \rho' \\ \sigma \times_{sh} \sigma' \end{pmatrix}, \begin{pmatrix} \rho'' \\ \sigma'' \end{pmatrix} \right\rangle \tag{7.12}$$

where on the left hand side the sum is over all good cuts $\sigma'' = \sigma_1'' * \sigma_2''$ and where it must be the case (by the definition of the multiplication) that $\mathrm{supp}(\rho) \cap \mathrm{supp}(\rho') = \emptyset$. Now on the right hand there can be a summand that is nonzero only if $\sigma'' = \rho * \rho'$ and then because $\rho$ and $\rho'$ have disjoint supports this is a good cut. Moreover $\rho''$ must be a shuffle of $\sigma, \sigma'$. If both these conditions hold the right hand side is 1, otherwise it is zero. Now on the left hand side there can be at most one good cut which yields a term that is nonzero, viz the one for which $\rho = \sigma_1''$ and $\rho' = \sigma_2''$. Further $\rho''$ is made up of the two complimentary subwords $p''^{-1}(\sigma_1'')$ and $p''^{-1}(\sigma_2'')$ with disjoint supports, that is, it is a shuffle of these two words. Now for the left hand side to have a term equal to one (there can be only one at most) it must first of all be the case that $\rho = \sigma_1''$ and $\rho' = \sigma_2''$ for some good cut $\sigma'' = \sigma_1'' * \sigma_2''$ so that also $\sigma'' = \sigma_1'' * \sigma_2''$ (and there can be at most one good cut like that). It must also be the case that $\sigma = p''^{-1}(\sigma_1'')$ and $\sigma' = p''^{-1}(\sigma_2'')$ making $\rho''$ a shuffle of $\sigma, \sigma'$. Thus if the left hand side of (7.10) is nonzero it is equal to one and then the right hand side is also equal to 1.

Inversely let the righthand side be equal to 1. Then $\sigma'' = \rho * \rho'$ and $\rho''$ must be a shuffle of $\sigma, \sigma'$. Let $\sigma_1''$ be the maximal subword of $\sigma''$ with the same support as $\sigma$ as a subword of $\rho''$. Then $\mathrm{supp}(\sigma_1'') = \mathrm{supp}(\sigma) = \mathrm{supp}(\rho)$ and as $\sigma'' = \rho * \rho'$ and $\mathrm{supp}(\rho) = \mathrm{supp}(\rho') = \emptyset$ it follows that $\sigma_1'' = \rho$ is actually a prefix of $\sigma''$ and that the left hand side of (7.10) is also equal to one.

There is a second structure on the Abelian group with as basis all substitutions *p*. This one is obtained by sort of imitating $Shuffle \otimes LieHopf$ (instead of $LieHopf \otimes Shuffle$ as was the case for the structure described above in this section.) It is defined as follows

$$m'(\begin{pmatrix} \rho \\ \sigma \end{pmatrix} \otimes \begin{pmatrix} \rho' \\ \sigma' \end{pmatrix}) = \begin{pmatrix} \rho \times_{sh} \rho' \\ \sigma * \sigma' \end{pmatrix} \tag{7.13}$$

$$\mu'(\begin{pmatrix} \rho \\ \sigma \end{pmatrix}) = \sum \begin{pmatrix} \rho_1 \\ p(\rho_1) \end{pmatrix} \otimes \begin{pmatrix} \rho_2 \\ p(\rho_2) \end{pmatrix} \tag{7.14}$$

where the sum in (7.14) is over all good cuts $\rho = \rho_1 * \rho_2$ and $p(\rho_i)$ denotes the maximal subword of $\sigma$ with the same support as $\rho_i$.

This Hopf algebra structure is isomorphic to the previous one. The isomorphism is given by

$$\begin{pmatrix} \rho \\ \sigma \end{pmatrix} \mapsto \begin{pmatrix} \sigma \\ \rho \end{pmatrix} \tag{7.15}$$

Define a (new) positive definite inner product structure $\langle \ , \ \rangle'$ on *dWHA* by declaring the substitutions to be an orthonormal basis. Then $m$ and $\mu'$ are dual to each other and so are $m'$ and $\mu$. That is

$$\langle p \otimes p', \mu(p'') \rangle' = \langle m'(p \otimes p'), p'' \rangle' \tag{7.14}$$

$$\langle p \otimes p', \mu'(p'') \rangle' = \langle m(p \otimes p'), p'' \rangle' \tag{7.15}$$

## 8. The word Hopf algebra *WHA*

This Hopf algebra is a sub Hopf algebra of *dWHA*. It has as basis substitutions of the following form

$$p = \begin{pmatrix} \rho \\ \sigma \end{pmatrix}, \quad \rho = [\underbrace{x_1, x_1, \cdots, x_1}_{r_1}, \underbrace{x_2, \cdots x_2}_{r_2}, \cdots, \underbrace{x_m, \cdots, x_m}_{r_m}] \tag{8.1}$$

That is besides the support condition on *p*, the top word has the property that if two letters are the same then all the letters between then are also equal to these. It is immediate to check that the product of such substitutions is again a sum of substitutions of this kind and also that the coproduct of a substitution of this form is a sum of tensor products of substitutions of this form. It will save typing (and printing ink and paper) to denote such a word $\rho$ as $[x_1^{r_1}, \cdots, x_m^{r_m}]$.

A substitution of the form (8.1) can be uniquely encoded as a single word over the integers as follows.

First let $\alpha = [a_1, a_2, \cdots, a_m]$ be a word over the natural numbers. Let

$$\{a_1', \cdots, a_n'\} = \mathrm{supp}(\alpha), \quad a_1' < \cdots < a_n'$$

Then the *WHA* substitution associated to $\alpha$ can be written

$$p(\alpha)=\begin{pmatrix}[x_{a'_1}^{r_1},\cdots,x_{a'_n}^{r_n}]\\ [x_{a_1},\cdots,x_{a_m}]\end{pmatrix} \tag{8.2}$$

where $r_1=a'_1,\cdots,r_i=a'_i-a'_{i-1},\cdots,r_n=a'_n-a'_{n-1}$. Inversely, if $p$ is a substitution of the form (8.1) then the word $\alpha$ associated to it is obtained as follows. Let $\sigma=[x_{i_1},\cdots,x_{i_t}]$, then $\alpha(p)=[a_1,\cdots,a_t]$ with $a_j=r_1+\cdots+r_{i_j}$.

For instance if $\alpha=[3,2,7,2,4]$

$$p(\alpha)=\begin{pmatrix}[x_2^2,x_3,x_4,x_7^3]\\ [x_3,x_2,x_7,x_2,x_4]\end{pmatrix}$$

and inversely, if $p$ is just like above, or written more canonically as

$$p=\begin{pmatrix}[x_1^2,x_2,x_3,x_4^3]\\ [x_2,x_1,x_4,x_1,x_3]\end{pmatrix}$$

$r_1=2,r_2=1,r_3=1,r_4=3$ and $i_1=2,i_2=1,i_3=4,i_4=1,i_5=3$ so that

$$\alpha(p)=[r_1+r_2,r_1,r_1+r_2+r_3+r_4,r_1,r_1+r_2+r_3]=[3,2,7,2,4]$$

The action of the basis elements of $WHA$ on $Shuffle$ is different from that mentioned in section 5. For example for the $p$ and $\alpha$ under discussion, $\mathrm{ht}(\alpha)=7$, it acts as zero on all words of length unequal to 7, and on the special words of length 7 of the form $[x_1,x_1,x_2,x_3,x_4,x_4,x_4]$ it acts exactly like in section 5, i.e. by picking out respectively the third, second, seventh, second, and fourth letter. The difference is that under this interpretation $\alpha$ is also zero on words of *Shuffle* of length 7 that are not of the special form $[x_1,x_1,x_2,x_3,x_4,x_4,x_4]$.

It is also obvious that the multiplication of these special substitutions, when written as words over the integers is like the one in section 1, see (1.4). That is

$$m_{WHA}(\alpha\otimes\beta)=[a_1,\cdots,a_m]\times_{sh}[\mathrm{ht}(\alpha)+b_1,\cdots,\mathrm{ht}(\alpha)+b_n]$$

It is trickier to write down a formula for the comultiplication; but a recipe is of course implied by the remarks above. Here is an example with $p$ and $\alpha$ as above. By definition

$$\mu_{WHA}(p)=1\otimes p+\begin{pmatrix}[x_2]\\ [x_2]\end{pmatrix}\otimes\begin{pmatrix}[x_1^2,x_3,x_4^3]\\ [x_1,x_4,x_1,x_3]\end{pmatrix}+\begin{pmatrix}[x_1^2,x_2,x_4^3]\\ [x_2,x_1,x_4,x_1]\end{pmatrix}\otimes\begin{pmatrix}[x_3]\\ [x_3]\end{pmatrix}+p\otimes 1$$

which translated to the level of words over the integers yields, with $\alpha=[3,2,7,2,4]$

$$\mu_{WHA}(\alpha)=1\otimes\alpha+[1]\otimes[2,6,2,3]+[3,2,6,2]\otimes[1]+\alpha\otimes 1$$

It is perhaps worth noting that with the interpretation given in this section of words acting as endomorphisms of *Shuffle* the dodgy procedure of section 5 for defining a comultiplication actually works; as it does for $dWHA$.

## 9. The permutation Hopf algebra *MPR* as a sub Hopf algebra of *dWHA*.

Let $\tau$ be a permutation word of length $n$. The permutation words over **N** are precisely the words of height equal to length and no multiplicities. To $\tau$ associate a substitution as follows

$$\varphi : \tau = [t_1, t_2, \cdots, t_n] \mapsto p(\tau) = \begin{pmatrix} [x_1, x_2, \cdots, x_n] \\ [x_{t_1}, x_{t_2}, \cdots, x_{t_n}] \end{pmatrix} \in dWHA \tag{9.1}$$

It is easy to characterize the substitutions that arise this way. They are precisely the substitutions for which the top word and the bottom word have no multiplicities.

More generally if $\tau = [t_1, t_2, \cdots, t_n]$ is any word over the natural numbers with no multiplicities there is also a permutation substitution attached to it. Indeed, let the support of such a $\tau$ be $\{a_1, a_2, \cdots, a_n\}$, $a_1 < a_2 < \cdots < a_n$ (same $n$ because of the no multiplicities condition). Then the permutation substitution associated to $\tau$ is

$$\begin{pmatrix} [x_{a_1}, x_{a_2}, \cdots, x_{a_n}] \\ [x_{t_1}, x_{t_2}, \cdots, x_{t_n}] \end{pmatrix}$$

which, as a substitution, is the same as $p(\mathrm{st}(\tau))$ as defined by (9.1), where st is the standardization map of section 1, i.e. $\mathrm{st}(\tau) = [\psi(t_1), \psi(t_2), \cdots, \psi(t_n)]$ where $\psi$ is the unique strictly monotone map $s\mathrm{upp}(\tau) \longrightarrow \{1, 2, \cdots, n\}$. This is how standardization of permutation words appears in this business. Of course also in the world of permutations it is quite customary to identify permutations defined by $\tau$ and $\mathrm{st}(\tau)$. (The underlying alphabeth does not really matter.)

Actually $p(\tau)$ is in $WHA \subset dWHA$ but at this stage of the investigations it seems better to regard $MPR$ primarily as a sub Hopf algebra of $dWHA$. This begs a quesion which is taken care of by the following proposition.

9.2. *Proposition.* The imbedding defined by (9.1) is a homogeneous monomorphism of connected graded Hopf algebras $MPR \longrightarrow dWHA$.

Here $MPR$ has the graded Hopf algebra structure defined in section 1, see especially (1.2), (1.3), and $dWHA$ has the Hopf algebra structure (the first one) described in section 7, see especially formulas (7.3)-(7.4).

*Proof.* That things go well for the grading and for the units and counits is immediate. So let

$$\sigma = [s_1, \cdots, s_m] \quad \text{and} \quad \tau = [t_1, \cdots, t_n]$$

and consider their associated substitutions according to (9.1). The product in $dWHA$ of these substitutions is

$$\begin{pmatrix} [x_1, \cdots, x_m] * [y_1, \cdots, y_n] \\ [x_{s_1}, \cdots, x_{s_m}] \times_{sh} [y_{t_1}, \cdots, y_{t_n}] \end{pmatrix} = \begin{pmatrix} [x_1, \cdots, x_m, x_{m+1}, \cdots, x_{m+n}] \\ [x_{s_1}, \cdots, x_{s_m}] \times_{sh} [x_{m+t_1}, \cdots, x_{m+t_n}] \end{pmatrix}$$

which is the sum of the permutation substitutions corresponding to $\sigma \times_{sh} [m + t_1, \cdots, m + t_n]$. This shows that $\varphi$ preserves multiplication. Now lets look at the comultiplication (in $dWHA$) on $p(\tau)$. This gives

$$\mu(p(\tau)) = \sum_{\text{good cuts}} \begin{pmatrix} p(\tau)^{-1}(\sigma_1) \\ \sigma_1 \end{pmatrix} \otimes \begin{pmatrix} p(\tau)^{-1}(\sigma_2) \\ \sigma_2 \end{pmatrix}, \quad p(\tau) = \begin{pmatrix} \rho \\ \sigma \end{pmatrix} \tag{9.3}$$

Now note that because there are no multiplicities all cuts are good cuts, and, again because there are no multiplicities, $\sigma_1$ is a permutation word on the alphabeth formed by the letters in $p(\tau)^{-1}(\sigma_1)$ which is a subword of $[x_1,\cdots,x_n]$ so that the $x$'s in $p(\tau)^{-1}(\sigma_1)$ appear in their natural order. Thus the sum of tensor products of permutation substitutions (9.3) corresponds to

$$\sum_{i=0}^{n} \mathrm{st}([t_1,\cdots,t_i]) \otimes \mathrm{st}([t_{i+1},\cdots t_n])$$

proving that $\varphi$ also preserves the comultiplications.

There is a second Hopf algebra structure on *MPR* defined as follows. Let $\sigma = [s_1,\cdots,s_m]$ and $\tau = [t_1,\cdots,t_n]$ be two permutation words, then

$$m'(\sigma \otimes \tau) = \sum u * v \tag{9.4}$$

where the sum is over all pairs of words $(u,v)$ such that $\mathrm{supp}(u) \cup \mathrm{supp}(v) = \{1,2,\cdots,m+n\}$, $\mathrm{st}(u) = \sigma$, $\mathrm{st}(v) = \tau$. And the comultiplication is given by

$$\mu'(\tau) = \sum_{i=0}^{n} \tau_{\{1,\cdots,i\}} \otimes \mathrm{st}(\tau_{\{i+1,\cdots,n\}})$$

where $\tau_I$ for $I \subset \{1,\cdots,n\}$ is the word obtained from the permutation word $\tau$ by retaining only the letters (digits) in *I*. E.g. if $\tau = [4,1,5,3,2]$, $\tau_{\{1,2\}} = [1,2]$, $\tau_{\{1,2,3\}} = [1,3,2]$, $\tau_{\{3,4,5\}} = [4,5,3]$, $\mathrm{st}(\tau_{\{3,4,5\}}) = [2,3,1]$.

The counit and unit are the same as for $(MPR, m, \mu, e, \varepsilon)$.

9.5. *Proposition.* The imbedding $\varphi$ given by (9.1) is (also) a homogeneous monomorphism of Hopf algebras $(MPR, m', \mu') \longrightarrow (dWHA, m', \mu')$. The isomorphism (7.13) $(dWHA, m, \mu) \longrightarrow (dWHA, m', \mu')$ induces an isomorphism $(MPR, m, \mu) \longrightarrow (MPR, m', \mu')$ given by $\tau \mapsto \tau^{-1}$ (as permutations).

*Proof.* It is not difficult to prove the first statement directly along similar lines as the proof of proposition 9.2. It is still easier to prove the second statement first and to use the known fact that $\tau \mapsto \tau^{-1}$ is an isomorphism between the two Hopf structures on *MPR*, [10], [11].

So let $\tau = [t_1,\cdots,t_n]$ be a permutation word. The corresponding substitution is (see (9.1))

$$p(\tau) = \begin{pmatrix} [x_1, x_2, \cdots, x_n] \\ [x_{t_1}, x_{t_2}, \cdots, x_{t_n}] \end{pmatrix} \in dWHA \tag{9.6}$$

Applying the isomorphism (7.13) to this one obtains

$$\begin{pmatrix} [x_{t_1}, x_{t_2}, \cdots, x_{t_n}] \\ [x_1, x_2, \cdots, x_n] \end{pmatrix}$$

As there are no multiplicities this is a permutation substitution. To find out to which permutation word it corresponds one has to relabel the entries of the top word so that the $x$'s are in the right order. This gives

$$\begin{pmatrix} [x_1, x_2, \cdots, x_n] \\ [x_{\tau^{-1}(1)}, x_{\tau^{-1}(2)}, \cdots, x_{\tau^{-1}(n)}] \end{pmatrix}$$

where $\tau^{-1}$ is the inverse permutation to $\tau$, i.e. $\tau^{-1}: t_i \mapsto i$. Now by loc. cit. the first statement of the proposition follows immediately.

Finally, there is the nondefinite inner product on $MPR$ defined by

$$\langle \sigma, \tau \rangle = \begin{cases} 1 & \text{if } \tau = \sigma^{-1} \\ 0 & \text{otherwise} \end{cases} \tag{9.7}$$

By what has just been said this is the restriction of the inner product (7.9). It follows that

9.8.*Proposition.* The Hopf algebras $(MPR, m, \mu)$ and $(MPR, m', \mu')$ are selfdual with respect to the inner product (9.7). Under the inner product $\langle\ ,\ \rangle'$ which makes the permutation words an orthonormal basis the Hopf algebra $(MPR, m, \mu)$ is dual to $(MPR, m', \mu')$.

This is just another way of getting a known result from [10], [11]. Note also that this way it is quite easy to prove that $(MPR, m, \mu)$ and $(MPR, m', \mu')$ are indeed Hopf algebras.

**10. More sub Hopf algebras of *dWHA* and *WHA*.**
There are quite a number of natural sub Hopf algebras of $WHA$ and $dWHA$. There are also still larger Hopf algebras of 'double word type'. Here are some of them.

*Injective words.* A word $\sigma = [s_1, s_2, \cdots, s_m]$ over the integers is called injective if their are no repeats (= no multiplicities), i.e if the cardinality of the support of $\sigma$ is equal to the length of $\sigma$ (which is $m$ in this case). For the associated substitution in $dWHA$

$$p = \begin{pmatrix} \rho \\ \sigma \end{pmatrix}$$

this is the same as saying that the bottom word of the substitution $p$ is multiplicity free. Call these injective substitutions. NB, this does not mean that the corresponding endomorphisms of *Shuffle* are injective; these endomorphisms are never injective (nor surjective for that matter). Obviously the product of two injective substitutions is a sum of injective substitutions. Also the coproduct of an injective substitution is obviousy a sum of tensor products of injective substitutions.

Thus the Abelian subgroup spanned by the injective words is a sub Hopf algebra of $WHA \subset dWHA$ which will be denoted $WHA_{inj}$.

*Surjective words.* A word $\sigma = [s_1, s_2, \cdots, s_m]$ is called surjective if there are no gaps. That means that the number of different natural numbers occurring in $\sigma$, the content $\mathrm{ct}(\sigma)$ of $\sigma$, is equal to the height of $\sigma$, which is the largest natural number occurring in $\sigma$. According to the recipe given in section 8 this means precisely that the associated substitution

$$p = \begin{pmatrix} \rho \\ \sigma \end{pmatrix}$$

has its top word multiplicity free. Call these surjective substitutions. Again it is virtually immediate, just looking at the defining formulas, that the product of two surjective substitutions is a sum of surjective substitutions and that the coproduct of a surjective substitutions is a sum of tensor products of surjective substitutions. Thus the surjective words span a sub Hopf algebra of $WHA \subset dWHA$ which will be denoted $WHA_{surj}$.

Note that $WHA_{inj} \cap WHA_{surj} = MPR$.

*Multisupport.* The multisupport of a word $\sigma = [s_1, s_2, \cdots, s_m]$ is the multiset of the different letters occurring in it together with their multiplicities. So technically it is a a set $S$ together with a function to $\mathbf{N}$ giving for each element of $S$ a multiplicity. If $S = \{s_1, \cdots, s_m\}$ is the underlying set of a multiset with $r_i$ the multiplicity of $s_i$; a useful notation for this multiset is $\{s_1^{r_1}, \cdots, s_m^{r_m}\}$. Thus for instance $\mathrm{msupp}([6,5,7,2,5,6,1,1,5]) = \{1^2, 2, 5^3, 6^2, 7\}$.

Now consider substitutions with top word $\rho$ and bottom word $\sigma$ such that

$$\mathrm{msupp}(\rho) = \mathrm{msupp}(\sigma)$$

Again it is an easy exercise to see that these substitutions define a sub Hopf algebra of $dWHA$. It will be denoted $dWHA_{msupp}$.

*Bounded multiplicity.* Fix a natural number $b \in \mathbf{N}$. Now consider substitutions for which the multiplicity of each letter that occurs is $\leq b$. These also span a sub Hopf algebra, which could be denoted $dWHA(b)$. A variant is to consider those substitutions for which the multiplicity of each letter in each of the two words is precisely $b$.

There are also generalizations of $dWHA$. For instance consider pairs of words $p = \begin{pmatrix} \rho \\ \sigma \end{pmatrix}$ such that the support of the top word is included in the support of the bottom word (but not neccessarily equal). Then the multiplication and comultiplication formula for $dWHA$ still make sense provided in the comultiplication formula one defines $p^{-1}(\sigma_i)$ as the maximal subword of $\rho$ whose support is included in the support of $\sigma_i$, $i = 1,2$.

Dually there is a generalization of $(dWHA, m', \mu')$ with the reverse inclusion of support condition. But I know of no ‘endomorphism interpretation’ for these generalizations.

All these word Hopf algebras (and there are many more, involving conditions on heigth and content and all kinds of combination conditions) have various inclusion relations. There are also natural projections. Here are two.

*Projecting $dWHA$ onto $MPR$.* In $dWHA$ consider all substitutions $p = \begin{pmatrix} \rho \\ \sigma \end{pmatrix}$ such that at least one letter has multiplicity $> 1$ in $\rho$ or $\sigma$ or both. Let $J_{mult}$ be the Abelian subgroup spanned by these substitutions. The claim is that $J_{mult}$ is a Hopf ideal. That means that if $p \in J_{mult}$ and $p'$ is any substitution then their product (in any order) is in $J_{mult}$ and $\mu(p) \in J_{mult} \otimes dWHA + dWHA \otimes J_{mult}$. The first is obvious because if there is an element of multiplicity 2 or larger in $\rho$ then there is an (actually the same) element of multiplicity 2 or larger in the concatenations $\rho * \rho'$, $\rho' * \rho$, and if there is an element of multiplicity two or larger in $\sigma$ the same is true in any shuffle of $\sigma$ and $\sigma'$. As to the comultiplication, if there is an element of multiplicity two or larger in $\sigma$ the same is true for $\sigma_1$ or $\sigma_2$ (or both) because only good cuts $\sigma = \sigma_1 * \sigma_2$ are involved, and if there is an element of multiplicity two or larger in $\rho$ then the same is true for $p^{-1}(\sigma_1)$ or $p^{-1}(\sigma_2)$ because these are complementary subwords of $\rho$ with disjoint supports (again because there are only good cuts involved). Note that if the basis element $p$ of $dWHA$ is not in $J_{mult}$ it defines a permutation, i.e. it is in $MPR \subset dWHA$.

Now define

$$\psi : dWHA \longrightarrow MPR, \quad p \mapsto \begin{cases} p & \text{if} \quad p \in MPR \subset dWHA \\ 0 & \text{if} \quad p \in J_{mult} \end{cases} \tag{10.1}$$

Because $J_{mult}$ is a Hopf ideal this is a morphism of Hopf algebras. It is also the identity on the sub Hopf algebra $MPR$ of $dWHA$ and hence gives a Hopf retraction of the inclusion $MPR \subset dWHA$.

*Surjective standardization.* Let $\sigma = [s_1, s_2, \cdots, s_m]$ be a word over the natural numbers. Part of the standardization process of [13] turns it into a surjective word. This goes as follows. Let

$$\mathrm{std}_\sigma : \mathrm{supp}(\sigma) \longrightarrow \{1, \cdots, \mathrm{ct}(\sigma)\}$$

be the unique strictly montone map between these ordered subsets of $\mathbf{N}$. Then

$$\mathrm{std}_{surj}: WHA \longrightarrow WHA_{surj} \tag{10.2}$$

is defined by

$$\sigma = [s_1, s_2, \cdots, s_m] \mapsto [\mathrm{std}_\sigma(s_1), \cdots, \mathrm{std}_\sigma(s_m)]$$

10.3. *Proposition.* The morphism of Abelian groups (10.2) is a Hopf algebra retraction of the inclusion $WHA_{surj} \subset WHA$ and induces a Hopf algebra retraction of $MPR \subset WHA_{inj}$.

Note that this retraction is different from the one defined by (10.1) combined with the inclusions $WHA_{inj} \subset WHA \subset dWHA$.

*Proof.* Routine, once it has been observed that the substitution $\begin{pmatrix} \rho' \\ \sigma' \end{pmatrix}$ defined by a word $\sigma$ has as the indices of the bottom word precisely the standardization of that word and that the substitution of the surjective standardization of $\sigma$ has the same bottom word and as top word the word obtained from $\rho'$ obtained by removing all multiplicities.

Note also that surjective standardization applied to (the indices of) both words of a substitution gives the same substitution.

*Standardization.* There is a process, due to Schensted, [13], that takes any word over the natural numbers into a permutation word. It goes as follows. Let $\sigma = [s_1, s_2, \cdots, s_m]$ be a word over $\mathbf{N}$. Let its multisupport be

$$\mathrm{msupp}(\sigma) = \{t_1^{r_1}, \cdots, t_n^{r_n}\}, \quad t_1 < \cdots < t_n$$

Now for each $i$ replace the $r_i$ entries of $\sigma$ that are equal to $t_i$ by the integers $r_1 + \cdots + r_{i-1} + 1, \cdots, r_1 + \cdots + r_i$ in this order, where $r_0 = 0$. There results a permutation word $\mathrm{st}(\sigma)$. For instance st([4,3,3,7,4,8,4])=[3,1,2,6,4,7,5]. Let $\mathrm{st}: WHA \longrightarrow MPR$ be the corresponding morphism of free Abelian groups. Then st is a morphism of algebras but not a morphism of coalgebras. So it defines an algebra retraction of the inclusion $MPR \subset WHA$ but not a Hopf algebra retraction.

Much more investigation of the relations between the various word Hopf algebras is called for. At the moment this is still in the infancy stage (or perhaps the foetal stage).

*Another grading.* For the double word Hopf algebra both the top and bottom word seem equally important. For all the others the bottom word seems dominant and one can define another grading by taking the length of the bottom word as degree. The result is that one obtains a connected graded Hopf algebra with the desirable property that all homogeneous components are of finite rank.

**11. Composition and cocomposition; second multiplications and comultiplications on word Hopf algebras.**

All of the sub Hopf algebras of *dWHA* are modules of endomorphism of *Shuffle* (or dually of *LieHopf*) and endomorphisms can be composed defining a second multiplication on these word Hopf algebras. As a rule these have no special distributivity properties; for instance they need not be distributive over the first multiplication in the Hopf algebra sense. For distributivity on the left of a Hopf algebra $(H, m_\Sigma, \mu_\Sigma)$ with second comultiplication $m_\Pi$ the following diagram needs to be commutative.

$$\begin{array}{ccc} H^{\otimes 3} & \xrightarrow{id\otimes m_\Sigma} & H^{\otimes 2} \\ \downarrow \mu_\Sigma \otimes \mathrm{id} \otimes \mathrm{id} & & \downarrow m_\Pi \\ H^{\otimes 4} & & H \\ \downarrow \mathrm{id}\otimes \mathrm{tw} \otimes \mathrm{id} & & \uparrow m_\Sigma \\ H^{\otimes 4} & \xrightarrow{m_\Pi \otimes m_\Pi} & H^{\otimes 2} \end{array} \tag{11.1}$$

Still it happens that composition induces a second multiplication that is distributive over the first. This happens for instance for the sub Hopf algebra of noncommutative symmetric functions *NSymm* of $MPR \subset dWHA$. See below for how *NSymm* is imbedded in *MPR*. In this case composition (on *MPR*, i.e.composition of permutations) induces a second multiplication on *NSymm* that is left distributive over the first one (but not right distributive). In a (at least at first sight) different way this also happens for the second multiplication on the functors of Witt vectors and generalized Witt vectors, see [7], section 12.

However, even if there are no distributive (or other nice properties) it can be most useful to have a second way of producing a new element for two given ones. See e.g. the use of such things in [6].

11.2. *Open problem.* What properties must there be on an element of a Hopf algebra with second multiplication in order that second mutiplication by that element is left distributive (right distributive) over the first multiplication in the Hopf algebra sense.

Here is how composition works out for *dWHA* and various sub Hopf algebras. Let

$$p = \begin{pmatrix} \rho \\ \sigma \end{pmatrix} \quad \text{and} \quad p' = \begin{pmatrix} \rho' \\ \sigma' \end{pmatrix}$$

be two substitutions. Then there product is zero unless there is a bijection $\varphi : \mathrm{supp}(\rho) \longrightarrow \mathrm{supp}(\sigma')$ taking $\rho$ into $\sigma'$. And in that case the the composition (second product) is

$$m_\Pi(p \otimes p') = p \circ p' = \begin{pmatrix} \rho' \\ \varphi(\sigma) \end{pmatrix} \tag{11.3}$$

where $\varphi(\sigma)$ is the word obtained by applying the change of variables $\varphi$ to the word $\sigma$

which has the same support as $\rho$, and where the interpretation as an endomorphism of *Shuffle* is: first $p'$ followed by $p$ (as indicated in (11.3).

It is easy to check that all the sub Hopf algebras $WHA$, $WHA_{inj}$, $WHA_{surj}$, $MPR$, $dWHA(b)$, $dWHA_{msupp}$ are stable under this second multiplication.

Dually there is comultiplication. Using the autoduality of $dWHA$ (or directly) this ought to be given by the formula

$$\mu_{\Pi}(\begin{pmatrix} \rho \\ \sigma \end{pmatrix}) = \sum \begin{pmatrix} \tau \\ \sigma \end{pmatrix} \otimes \begin{pmatrix} \rho \\ \tau \end{pmatrix} \tag{11.4}$$

where the sum is over all $\tau$ with the same support as $\rho$ and $\sigma$ and such that each of the two factors in (11.4) are in the appropriate sub Hopf algebra. However in the case of $dWHA$, and $WHA$ there are infinitely many such $\tau$ and comultiplication is not well defined. In the case of the sub Hopf algebras $WHA_{inj}$, $WHA_{surj}$, $MPR$, $dWHA(b)$, $dWHA_{msupp}$ there are only finitely may such $\tau$ and the second comultiplication is well defined.

**12.** *Symm, NSymm, QSymm,* **and** *MPR*: **part of the diagram of their relations.**

Most of the rest of this paper is about the Hopf algebras occurring in the following diagram and their relations.

$$\begin{array}{ccccc} (MPR, m', \mu') & & & Symm \\ i \uparrow & \ddots & \cdot^{\cdot^{\cdot}} & \downarrow i' \\ NSymm & \cdots & \cdots & QSymm \\ \pi' \downarrow & \cdot^{\cdot^{\cdot}} & \ddots & \uparrow \pi \\ Symm & & & (MPR, m, \mu) \end{array} \tag{12.1}$$

where two $\cdots$'s, either horizontal or diagonal indicate a duality pairing. The Hopf algebra morphisms $\pi, \pi'$ are surjective and the Hopf algebra morphisms $i, i'$ are injective; $\pi$ and $i$ are dual to each other and so are $\pi'$ and $i'$. The part of the diagram not involving $MPR$ has been discussed in some detail in [7, 8]. So here most of the emphasis will be on the relations involving one of the two isomorphic versions of $MPR$. The first one of these has been described in section 1; the second one in section 9.

One can of course change the diagonal duality $\langle \, , \, \rangle'$ between $(MPR, m', \mu')$ and $(MPR, m, \mu)$ to the autoduality defind by the non positive definite pairing of $(MPR, m, \mu)$ given by (9.7). A similar diagram ensues with $i$ replaced by $(\sigma \mapsto \sigma^{-1}) \circ i$ and $(MPR, m', \mu')$ in the upper left hand corner replaced by $(MPR, m, \mu)$ (or, equivalently, replacing $\pi$ by $\pi \circ (\sigma \mapsto \sigma^{-1})$ and $(MPR, m, \mu)$ by $(MPR, m', \mu')$ in the lower right hand corner).

For completeness sake (and the convenience of the reader) also $NSymm$ and $QSymm$ will be very briefly described below.

*NSymm*. As an algebra the Hopf algebra of noncommutative symmetric functions $NSymm$ is simply the free associative algebra in countably many (noncommuting) indeterminates

$$NSymm = \mathbf{Z}\langle Z_1, Z_2, \cdots \rangle \tag{12.2}$$

The comultiplication is determined by

$$Z_n \mapsto \sum_{i+j=n} Z_i \otimes Z_j \quad \text{where} \quad Z_0 = 1 \tag{12.3}$$

(and requiring this to be an algebra morphism). This is welldefined because *NSymm* is free as an algebra. So there is no difficulty establishing the Hopf property. The unit is 1 and the counit takes 1 to 1 and $\varepsilon(Z_n) = 0,\ n = 1,2,\cdots$ .

One basis of the underlying free Abelian group is formed by the noncommutative monomials

$$Z_\alpha = Z_{a_1} Z_{a_2} \cdots Z_{a_m}, \quad \alpha = [a_1, a_2, \cdots a_m] \in \mathbf{N}^*$$

indexed by all words over the natural numbers $\mathbf{N} = \{1,2,\cdots\}$. In this context these words are often called compositions, more precisely $\alpha$ is called a composition of $n = \mathrm{wt}(\alpha) = a_1 + \cdots + a_m$. *NSymm* is graded by giving $Z_\alpha$ degree $\mathrm{wt}(\alpha)$. Two other bases will play a role later on. Let the $S_n$ be defined by the noncommutative Wronski relations

$$S_n - S_{n-1}Z_1 + S_{n-2}Z_2 + \cdots + (-1)^{n-1} S_1 Z_{n-1} + (-1)^n Z_n = 0 \tag{12.4}$$

(or the same relations with the $Z$'s and $S$'s reversed; it does not matter which is used: the resulting elements are the same). If the $Z$'s are viewed as a noncommutative analogue of the complete symmetric functions, the $S$'s are the corresponding noncommutative elementary symmetric functions (or vice versa). Obviously the noncommutative monomials in the $S$'s also form a basis for the underlying Abelian group of *NSymm*. A third basis is given by the socalled ribbon Schur functions

$$R_\alpha = \sum_{\alpha \geq \beta} (-1)^{\lg(\alpha) - \lg(\beta)} S_\beta \tag{12.5}$$

where $\alpha \geq \beta$ means that $\alpha = [a_1, \cdots, a_m]$ is a refinement of $\beta = [b_1, \cdots b_n]$ (or $\beta$ a coarsening of $\alpha$)[8]. This is the case if there are indices $1 \leq j_1 < j_2 < \cdots < j_n = m$ such that $b_i = a_{j_{i-1}+1} + \cdots + a_{j_i}$ , which can be pictorially represented as

$$[\underbrace{a_1, \cdots, a_{j_1}}_{b_1}, \underbrace{a_{j_1+1}, \cdots, a_{j_2}}_{b_2}, \cdots, \underbrace{a_{j_{m-1}+1}, \cdots, a_m}_{b_n}]$$

*QSymm*. The Hopf algebra of quasisymmetric functions *QSymm* is the graded dual of *NSymm*. A basis for the underlying Abelian group is formed by all words $\alpha$, $\alpha \in \mathbf{N}^*$ with duality pairing

$$\langle Z_\beta, \alpha \rangle = \delta_\beta^\alpha \quad \text{(Kronecker delta)} \tag{12.6}$$

The mutiplication is the overlapping shuffle product, which can be described as follows. Take two words $\alpha = [a_1, \cdots, a_m]$ and $\beta = [b_1, \cdots b_n]$. For each $k \leq \min\{m, n\}$ take a 'word' with $m + n - k$ so far empty slots. Choose $m$ slots and in these put the entries of $\alpha$ in their original order; choose $k$ of these now filled slots; together with the $n - k$ still empty slots this makes $n$ slots; in these slots place the entries of $\beta$ in their original order; finally, for those slots that have two entries add them. The overlapping shuffle product, $\alpha \times_{osh} \beta$ , of $\alpha$ and $\beta$ is the sum of all words (compositions) that can be obtained this way. For instance

[8] The opposite sign is also used in the literature. I prefer this one and use as a mnemonic that if one composition refines another it has more parts.

$$[2]\times_{osh}[1,5,3]=[2,1,5,3]+[1,2,5,3]+[1,5,2,3]+1,5,3,2]$$
$$+[3,5,3]+[1,7,3]+[1,5,5]$$

$$[2]\times_{osh}[2,2,3]=3[2,2,2,3]+[2,2,3,2]+[4,2,3]+[2,4,3]+[2,2,5]$$

$$[1,1]\times_{osh}[1,1]=6[1,1,1,1]+2[1,1,2]+2[1,2,1]+2[2,1,1]+[2,2]$$

The empty word is the unit element.
The comultiplication is 'cut':

$$\mu([a_1,\cdots,a_m]=\sum_{i=0}^{m}[a_1,\cdots,a_i]\otimes[a_{i+1},\cdots,a_m] \tag{12.7}$$

and the counit is 1 on the empty composition and zero on all others.
Another basis of the underlying Abelian group is formed by the sums

$$F_\alpha=\sum_{\beta\geq\alpha}\beta \tag{12.8}$$

For instance $F_{[1,1,1]}=[1,1,1]$, $F_{[3]}=[3]+[2,1]+[1,2]+[1,1,1]$, $F_{[1,2]}=[1,2]+[1,1,1]$. These $F$'s are biorthogonal to the ribbon Schur functions under the duality pairing

$$\langle R_\alpha,F_\beta\rangle=\delta_{\alpha,\beta}\quad\text{(Kronecker delta)} \tag{12.9}$$

(This is not difficult; or see [5].)

*Compositions and descent sets*. Let $\sigma=[s_1,s_2,\cdots,s_m]$ be a permutation (word). Its descent set, $\text{desc}(\sigma)$, is the subset of $\{1,2,\cdots,m-1\}$ consisting of those natural numbers $i$ for which $s_i>s_{i+1}$. For instance $\text{desc}([3,2,5,7,1,4,6])=\{1,4\}\subset\{1,2,\cdots,6\}$. It is important to keep in mind that a descent set is not just a set, but a subset of a specific $\{1,2,\cdots,n\}$; thus $\{1\}\subset\{1,2\}$ is a very different object from $\{1\}\subset\{1,2,3,4\}$.

Given a descent set $D=\{d_1<\cdots<d_r\}\subset\{1,2,\cdots,m-1\}$ there is an associated composition, $\text{comp}(D)=[d_1,d_2-d_1,\cdots,d_r-d_{r-1},m-d_r]$ of $m$. Inversely, given a composition $\alpha=[a_1,a_2,\cdots,a_r]$ of $m$ there is an associated descent set, $\text{desc}(\alpha)=\{a_1,a_1+a_2,\cdots,a_1+\cdots+a_{r-1}\}\subset\{1,\cdots,m-1\}$. Note that $\text{desc}(\sigma)$ for $\sigma$ a permutation and $\text{desc}(\alpha)$ for $\alpha$ a composition are very different things.

The two recipes that associate a descent set to a composition and a composition to a descent set are inverse to each other.

*The injection of Hopf algebras* $i:\ NSymm\longrightarrow(MPR,m',\mu')$. Consider the morphism of algebras determined by

$$S_n\mapsto[1,2,\cdots,n] \tag{12.10}$$

i.e. $S_n$ goes to the identity permutation on $n$ letters. Because the $S_n$ obviously are also a free associative set of generators for $NSymm$ this is a well defined morphism of algebras. It is in fact a morphism of Hopf algebras. It sufices to check this on the generators $S_n$ (because of the Hopf property). Now

$$\mu_{NSymm}(S_n)=\sum_{i+j=n}S_i\otimes S_j,\quad S_0=1$$

as follows easily from the Wronski relations (12.4). On the other hand

$$\mu'_{MPR}([1,2,\cdots,n]) = \sum_{i=0}^{n}[1,2,\cdots,n]_{\{1,\cdots,i\}} \otimes \mathrm{st}([1,2,\cdots,n]_{\{i+1,\cdots,n\}})$$
$$= \sum_{i=0}^{n}[1,2,\cdots,i] \otimes [1,2,\cdots,n-i]$$

proving that this a Hopf algebra morphism.

To see that this morphism is injective a little more is needed. Recall that the multiplication $m'_{MPR}$ is given by the formula

$$m'_{MPR}(\rho \otimes \sigma) = \sum u * v$$

where the sum is over all concatenations $u * v$ such that

$$\mathrm{st}(u) = \rho = [r_1,\cdots,r_m],\ \mathrm{st}(v) = \sigma = [s_1,\cdots,s_n]$$

and

$$\mathrm{supp}(u) \cup \mathrm{supp}(v) = \{1,2,\cdots,n+m\},$$

where $m = \lg(\rho)$, $n = \lg(\sigma)$. It follows immediately that the descent set of any of the summands is either $\mathrm{desc}(\rho) \cup m + \mathrm{desc}(\sigma)$ or $\mathrm{desc}(\rho) \cup \{m\} \cup m + \mathrm{desc}(\sigma)$ depending on whether the last letter of $u$ is smaller or larger than the first letter of $v$. And in fact both possibilities always occur; the first one for $u = [r_1,\cdots,r_m]$, $v = [m+s_1,\cdots,m+s_n]$ and the second one for $u = [n+r_1,\cdots,n+r_m]$, $v = [s_1,\cdots,s_n]$.

It follows immediately that the largest (in size) descent set that occurs in the iterated product

$$m'_{MPR}([1,\cdots,a_1] \otimes [1,\cdots,a_2] \otimes \cdots \otimes [1,\cdots,a_m])$$

of identity permutations is $\{a_1, a_1 + a_2, \cdots, a_1 + \cdots + a_{m-1}\}$ and this one occurs for the following summand of this iterated product, $n = \mathrm{wt}(\alpha) = a_1 + \cdots + a_m$

$$[n-a_1+1, n-a_1+2,\cdots,n; n-a_1-a_2+1,\cdots,n-a_1;\cdots;$$
$$n-a_1-a_2-\cdots-a_m+1 = 1,\cdots,n-a_1-a_2-\cdots a_{m-1} = a_m]$$

These are all different and this proves that the morphism defined by (12.10) is injective.[9]

Let $D \subset \{1,\cdots,m-1\}$ be a subset. The descent class sum corresponding to it is

$$\vartheta_D = \sum_{\mathrm{desc}(\sigma)=D} \sigma \ \in\ MPR \qquad (12.11)$$

where the sum is over all permutations of $\{1,2,\cdots,m\}$ with descent set $D$; let $D' \subset \{1,2,\cdots,n\}$ be a second descent set, and let $D_1 = D \cup m + D' \subset (1,2,\cdots,m+n-1\}$, $D_2 = D \cup \{m\} \cup m + D'$. Then

---

[9] As a matter of fact these permutations for $m \geq 2$ are all part of a free generating set of $MPR$, see [11].

$$m'_{MPR}(\vartheta_D \otimes \vartheta_{D'}) = \vartheta_{D_1} + \vartheta_{D_2} \tag{12.12}$$

To see this, first note that by definition there are no multiplicities. Indeed if $uv = u'v'$ and $u$ and $u'$ are of equal length then $u = u'$ and $v = v'$. Thus, by the remarks above, the sum on the left hand side of (12.12) is part of the sum on the right hand side. Inversely, let $\tau$ be a permutation of $\{1,2,\cdots,m+n\}$ with descent set $D_1$ or $D_2$. Let $\tau = u * v$ be the concatenation decomposition (cut) of $\tau$ with $u$ of length $m$. Then $D = \mathrm{desc}(\mathrm{st}(u))$ and $D' = \mathrm{desc}(\mathrm{st}(v))$, so that $\tau$ also occurs on the left hand side of (12.12). This proves this formula.

The ribbon Schur functions in *NSymm* multiply as follows, [5],

$$R_{[a_1,\cdots,a_m]} R_{[b_1,\cdots,b_n]} = R_{[a_1,\cdots,a_m,b_1,\cdots,b_n]} + R_{[a_1,\cdots,a_{m-1},a_m+b_1,\cdots,b_n]} \tag{12.13}$$

Also note that if $D \subset \{1,2,\cdots,m-1\}$ is the descent set associated to a composition $[a_1,a_2,\cdots,a_r]$ of $m$ and $D' \subset \{1,2,\cdots,n-1\}$ is the descent set associated to a composition $[b_1,b_2,\cdots,b_n]$, then

$$\begin{aligned}
&\mathrm{desc}([a_1,a_2,\cdots,a_r,b_1,b_2,\cdots,b_s]) \\
&= \{a_1,a_1+a_2,\cdots,a_1+\cdots+a_{r-1},a_1+\cdots+a_r = m, m+b_1, m+b_1+b_2,\cdots,m+b_1+\cdots+b_{s-1}\} \\
&= D \cup \{m\} \cup m + D' \subset \{1,2,\cdots m+n-1\}
\end{aligned}$$

$$\begin{aligned}
&\mathrm{desc}([a_1,a_2,\cdots,a_{r-1},a_r+b_1,b_2,\cdots,b_s]) \\
&= \{a_1,a_1+a_2,\cdots,a_1+\cdots+a_{r-1},a_1+\cdots a_r+b_1, m+b_1+b_2,\cdots,m+b_1+\cdots b_{s-1}\} \\
&= D \cup m + D' \subset \{1,2,\cdots,m+n-1\}
\end{aligned}$$

Combining this with (12.12) and (12.13) leads to the following explicit description of the imbedding (12.10) of *NSymm* into *MPR*.

$$R_\alpha \mapsto \sum_{\mathrm{desc}(\alpha)=\mathrm{desc}(\sigma)} \sigma \tag{12.14}$$

Note that there are two different meanings of 'desc' involved here: $\mathrm{desc}(\sigma)$ is the descent set of a permutation while $\mathrm{desc}(\alpha)$ is the descent set corresponding to the composition $\alpha$.

A descent class in *MPR* is the set of all permutations (of the same length) with the same descent set. So the right hand side of (12.12) is a sum over a descent class. These form a basis for the image of (12.10), (12.12). This image, which, by what has been said, is isomorphic to *NSymm* as a Hopf algebra is known as the Solomon descent algebra. It carries a second multiplication given by composition of permutations that is left distributive over the first multiplication. The fact that the composition product of two sums over descent classes is a sum of sums over descent classes was discovered by Solomon, [15].

As the composition of two identity permutations is zero (if they are of uequal length) or that same identity permutation (if they are of equal length, it follows that the second multiplication on *NSymm* satisfies

$$\begin{aligned}
&m_\Pi(S_n \otimes S_n) = S_n \\
&m_\Pi(S_\alpha \otimes S_\beta) = 0 \ \text{ if } \ \mathrm{wt}(\alpha) \neq \mathrm{wt}(\beta)
\end{aligned} \tag{12.15}$$

By left distributivity and the grading properties this suffices to characterize the second multiplication.

*The projection* $\pi : (MPR, m, \mu) \longrightarrow QSymm$. The imbedding $NSymm \longrightarrow (MPR, m', \mu')$ by graded duality induces a projection $\pi : (MPR, m, \mu) \longrightarrow QSymm$. This one is explicitly given by

$$\sigma \mapsto F_{\mathrm{comp}(\mathrm{desc}(\sigma))} \tag{12.16}$$

This follows from (12.12) because the $F$'s and $R$'s are dual bases, see (12.9). The projection (12.16) turns out to be a very natural one[10]. The algebra of quasisymmetric functions is, by definition of 'quasisymmetric function', a subalgebra of the commutative power series over $\mathbf{Z}$ in variables $x_1, x_2, \cdots$ . There is an interpretation of $MPR$ as a subalgebra of the noncommuting powr series over the integers in these variables, and then (12.14) is induced by Abelianization, see [4].

**13. lsd permutations**.

It will be convenient to also consider permutations of other subsets of $\mathbf{N}$ then the $\{1,...,n\}$. Thus any injective word can be seen as defining a permutation of its support set. The descent set of such a permutation $\sigma = [b_1, b_2, \cdots, b_n]$ is defined as before $i \in \mathrm{desc}(\sigma)$ iff $b_i > b_{i+1}$. Let $D$ be any nonempty subset of $\{1,2,\cdots,n-1\}$ and write it in the form

$$D = \{i_1, i_1+1, \cdots, i_1+j_1-1; i_2, i_2+1, \cdots, i_2+j_2-1; \cdots; i_r, i_r+1, \cdots, i_r+j_r-1\} \tag{13.1}$$

with $i_s, j_s \in \mathbf{N}$, $i_s - (i_{s-1} + j_{s-1} - 1) \geq 2$ . For instance if $D = \{1,2,6,7\} \subset \{1,2,\cdots,7\}$ , $i_1 = 1, j_1 = 2, i_2 = 6, j_2 = 2, r = 2$ .

Now consider an alphabet $\{a_1 < a_2 < \cdots < a_{n-1} < a_n\} \subset \mathbf{N}$ of size $n$ and consider the following permutation of it

[10] See also section 17 below.

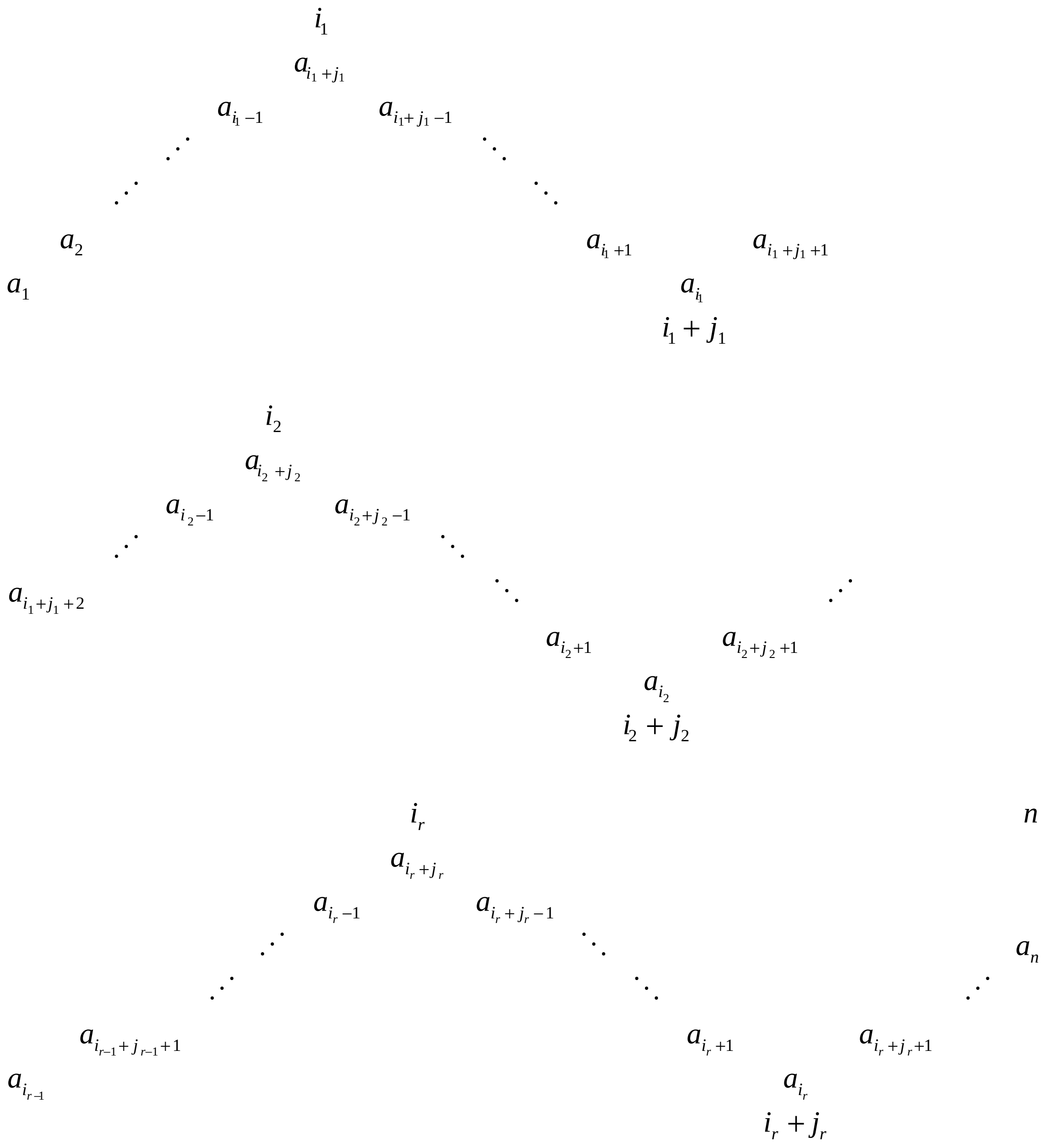

(13.2)

where the relative height positioning indicates whether one is at a point of ascent or descent and where the occasional numbers in the bottom and top lines indicate position; thus for example the letter $a_{i_1}$ is the $(i_1 + j_1) -$th letter of this permutation. Of course the initial 'up run' may be missing (if $i_1 = 1$) and the last up run may also be missing (if $i_r + j_r = n$).

The descent set of the permutation (13.2) is the set (13.1).

13.3. *Theorem.* The permutation (13.2) is the lexicographically smallest permutation in the descent class of all permutations of the alphabet $\{a_1, \cdots, a_n\}$ with descent set $D$ as given by (13.1).[11]

The lexicographically smallest permutation of a given descent class I call an lsd permutation (lsd is of course an acronym for **l**exicographically **s**mallest in its **d**escent class). If $D$ is empty (this corresponds to $r = 0$) there is only one permutation with descent class $D$. This is the identity permutation.

[11] It is also the smallest element in its descent class under the weak order, and each descent class has a unique smallest element under the weak order; see section 15 below.

*Proof.* This is not a difficult matter and handled by induction. The induction starts because things are obvious for an alphabet of size 1 or size 2.

There are two cases to consider:

A) $1 \in D$, i.e. $i_1 = 1$

B) $1 \notin D$, i.e. $i_1 > 1$

In the second case the permutation (13.2) starts with $a_1$ the smallest element of the alphabet involved. Removing this element from the alphabeth and from the permutation (13.2) we find a permutation of the same type for the descent set $D-1 \subset \{1,2,\cdots,n-2\}$. By induction (13.2) with this element removed is lexicographically smallest in its descent class and then, because $a_1$ is the smallest letter of the alphabet, (13.2) is the lexicographically smallest element defined by $D$.

In case A) the descent set starts with $j_1$ consecutive descents. So the first element of any permutation with this descent set must start with $a_{j_1+1}$ or a larger element. The permutation (13.2) in this case does indeed start with $a_{j_1+1}$. Remove this element from the alphabeth and from (13.2). There remains a permutation of the same shape with descent set $(D \setminus \{1\})-1 \subset \{1,2,\cdots,n-2\}$. With induction this one is lexicographically smallest in its descent set. It follows that (13.2) is lexicographically smallest in the descent set defined by $D$.

Another way to write the permutation (13.2) is as follows

$$\left(\begin{matrix} a_1 & a_2 & \cdots & a_{i_1-1} & ; & a_{i_1} & a_{i_1+1} & \cdots & a_{i_1+j_1} & ; & a_{i_1+j_1+1} & \cdots & a_{i_2-1} & ; \\ a_1 & a_2 & & a_{i_1-1} & ; & a_{i_1+j_1} & a_{i_1+j_1-1} & \cdots & a_{i_1} & ; & a_{i_1+j_1+1} & \cdots & a_{i_2-1} & ; \end{matrix}\right.$$

$$\begin{matrix} a_{i_2} & \cdots & a_{i_2+j_2-1} & a_{i_2+j_2} & ; & \cdots & ; & a_{i_{r-1}+j_{r-1}+1} & \cdots & a_{i_r-1} & ; \\ a_{i_2+j_2} & \cdots & a_{i_2+1} & a_{i_2} & ; & \cdots & ; & a_{i_{r-1}+j_{r-1}+1} & \cdots & a_{i_r-1} & ; \end{matrix} \tag{13.4}$$

$$\left.\begin{matrix} a_{i_r} & \cdots & a_{i_r+j_r} & ; & a_{i_r+j_r+1} & \cdots & a_n \\ a_{i_r+j_r} & \cdots & a_{i_r} & ; & a_{i_r+j_r+1} & \cdots & a_n \end{matrix}\right)$$

which makes it obvious that a non identity lsd permutation is always an involution. But not every involution is lsd.

13.4. *Corollary*. Let $u * v = \sigma$ be a cut of an lsd permutation. Then both $u$ and $v$ are lsd permutations (in their respective alphabets).

This is obvious from the form of (13.4). Note that if the cut takes place somewhere in an 'up run' then the alphabet of $u$ is an initial chunk from the alphabet $\{a_1,\cdots,a_n\}$, but if the cut occurs in a down run that is not the case. For example for the lsd permutation

$$\sigma = \begin{pmatrix} 1 & 2 & ; & 3 & 4 & 5 & 6 & ; & 7 & ; & 8 & 9 & 10 & ; & 11 & 12 \\ 1 & 2 & ; & 6 & 5 & 4 & 3 & ; & 7 & ; & 10 & 9 & 8 & ; & 12 & 11 \end{pmatrix}$$

and the cut (bottom line of $\sigma$ )$= u * v$ with $u$ of length 9 and $v$ of length 3, the alphabet of $u$ is $\{1,2,3,4,5,6,7,9,10\}$ and that of $v$ is $\{8,11,12\}$.

## 14. Retractions of $i : NSymm \longrightarrow MPR$ and sections of $\pi : MPR \longrightarrow QSymm$.

In this section the question is examined whether there are retractions of $i : NSymm \longrightarrow MPR$ that are algebra morrphisms and, dually, whether there are sections of $\pi : MPR \longrightarrow QSymm$ that are coalgebra morphisms. The first coalgebra section of $\pi : MPR \longrightarrow QSymm$ was constructed in [1]. Here several others will be described as well as an algebra retraction of $i : NSymm \longrightarrow MPR$.

Let $J_{nonlsd}$ be the subgroup of $MPR$ spanned by all those permutations that are non lsd.

14.1. *Proposition.* $J_{nonlsd}$ is an ideal in the algebra $(MPR, m')$.

*Proof.* This is an immediate consequence of corollary 13.4. Indeed let $\rho$ and $\sigma$ be two permutations and let at least one of them be non lsd. Consider a term $u * v$ from the product $m'(\rho \otimes \sigma)$. If it were lsd both $u$ and $v$ would be lsd. But a permutation $u$ is lsd if and only if st(*u)* is lsd. But $\mathrm{st}(u)$ must be equal to $\rho$ and $\mathrm{st}(v) = \sigma$, proving the proposition by contradiction.

14.2. *Theorem.* Define $\psi : (MPR, m') \longrightarrow NSymm$ by

$$\psi(\sigma) = \begin{cases} 0 \ \text{ if } \ \sigma \text{ is not an lsd permutation} \\ R_{\mathrm{comp}(\mathrm{desc}(\sigma))} \ \text{ if } \ \sigma \text{ is an lsd permutation} \end{cases} \tag{14.3}$$

Then $\psi$ is a morphism of algebras that is an algebra retraction of the imbedding of Hopf algebras $i : NSymm \longrightarrow (MPR, m')$.

*Proof.* By definition, see also (12.14), $\psi \circ i = \mathrm{id}_{NSymm}$. So it remains to show that $\psi$ is multiplication preserving.

First let $\sigma$ and $\tau$ be two permutations (basis elements of $MPR$) of which at least one is not lsd. Then their product is in the ideal $J_{nonlsd}$ on which $\psi$ is zero. So $\psi$ is multiplicative in this case.

Now let $\sigma$ and $\tau$ be two permutations that are both lsd. Let their two descent sets be $D \subset \{1,2,\cdots,m-1\}$` and $D' \subset \{1,2,\cdots,n-1\}$. As before, let $D_1 = D \cup m + D' \subset \{1,2,\cdots,m+n\}$ and $D_2 = D \cup \{m\} \cup m + D'$. Let $\rho_i$ be the unique lsd permutation with descent set $D_i$, $i = 1,2$. Then, taking into account the formulas for descent sets just below (12.13), it sufffices to show that both $\rho_i$ occur as summands in $m'(\sigma \otimes \tau)$ (also taking into account that there are no multiplicities, i.e that $m'(\sigma \otimes \tau)$ is a sum of different permutations with coefficient 1). So let $\rho_i = u_i * v_i$ be the cut with the length of $u_i$ equal to $m$. Then, by corollary 13.4 both $u_i$, $v_i$ are lsd and hence so are $\mathrm{st}(u_i)$, $\mathrm{st}(v_i)$ and these have descent sets $D$, $D'$ .Hence $\mathrm{st}(u_i) = \sigma$, $\mathrm{st}(v_i) = \tau$ proving what is needed.

It is also not difficult to write down the $\rho_i$ explicitly. Indeed let

$$\sigma = [a_1, \cdots, a_{r-1}, a_r = m, m-1, \cdots, r]$$

$$\tau = [s, s-1, \cdots, 1, b_{s+1}, \cdots, b_n]$$

(Because these permutations are lsd they must look like this, where of course it can be that $m = r$ (so that $\sigma$ ends with an up run) and it can be that $s = 1$ so that $\tau$ starts with an up run.) Then

$$\rho_1 = [a_1, \cdots, a_{r-1}, a_r = m, m-1, \cdots, r, m+s, m+s-1, \cdots, m+1, m+b_{s+1}, \cdots, m+b_n]$$

is lsd with descent set $D_1$ and

$$\rho_2 = [\underbrace{a_1, a_2, \cdots, a_{r-1}}_{\substack{\text{alphabet} \\ \{1,2,\cdots,r-1\}}}, \underbrace{m+s, m+s-1, \cdots, r+s, r+s-1, \cdots, r}_{\substack{\text{alphabet} \\ \{r,\cdots,m+s\}}}, \underbrace{m+b_{s+1}, \cdots, m+b_n}_{\substack{\text{alphabet} \\ \{m+s+1,\cdots,m+n\}}}]$$

is lsd with descent set $D_2$.

Let $MPR_{lsd}$ be the free Abelian subgroup spanned by the lsd permutations. Then by corollary 13.4 $MPR_{lsd}$ is a subcoalgebra of $(MPR, \mu)$[12].

14.4. *Theorem.* Consider the following commutative diagram

$$\begin{array}{ccc} MPR_{lsd} & \xrightarrow{\subset} & MPR \\ \downarrow \pi'' & & \downarrow \pi \\ QSymm & = & QSymm \end{array} \tag{14.5}$$

where $\pi''$ is defined as the restriction of $\pi$ to $MPR_{lsd}$. Then $\pi''$ is an isomorphism of coalgebras and its inverse is a coalgebra section of the projection of Hopf algebras $\pi: (MPR, m, \mu) \longrightarrow QSymm$.

*Proof.* As the composition of two coalgebra morphisms $\pi''$ is certainly a morphism of coalgebras. Moreover for every composition $\alpha$ of $m$ there is precisely one lsd permutation $\sigma$ of $\{1,2,\cdots,m\}$ such that $\mathrm{comp}(\mathrm{desc}(\sigma)) = \alpha$ so that $\pi''$ sets up a bijection between the basis of lsd permutations of $MPR_{lsd}$ and the $F$-basis of $QSymm$ making it an isomorphism.

The coalgebra section of theorem 14.4 is different from the coalgebra section which arises by duality from theorem 14.2 and both are different from the one in [1].[13]

Under suitable circumstances a surjective morphism of Hopf algebras $H \longrightarrow K$ which has a coalgebra section gives rise to a description of $H$ as a crossed product of $K$ with a certain kernel, see [2], see also [3].

14.6. *Open problem.* Prove a crossed product type theorem (BCM theorem) over the integers for connected graded Hopf algebras. Describe the situation for an injection of such Hopf algebras which admits an algebra retraction.

## 15. Combinatorics of the weak order and descent sets.

Let $\mathrm{Inv}(\sigma) = \{(i, j) \in \{1,\cdots,n\} \times \{1,\cdots,n\} : i < j \text{ and } s_i > s_j\}$ be the set of inversions of a permutation (word) $\sigma = [s_1, \cdots, s_n] \in S_n$. The left weak order on $S_n$ is defined by

$$\sigma <_{lwo} \tau \iff \mathrm{Inv}(\sigma) \subset \mathrm{Inv}(\tau) \tag{15.1}$$

This defines a partial order on $S_n$ which has a smallest element, viz the identity permutation (with empty set of inversions) and a largest element, viz $[n, n-1, \cdots, 2, 1]$.

For any partially ordered set $P$ its Hasse diagram is the graph with vertices the elements of $P$ drawn such that an element $u$ is lower than $v$ if and only if $u$ is smaller than $v$ and with an edge from $v$ to $u$ iff $u$ is smaller than $v$ and there are no elements between $v$ and $u$. This is called the covering relation of the partial order.

[12] It is also a subcoalgebra of $(MPR, \mu')$ because the isomorphism $(MPR, \mu) \longrightarrow (MPR, \mu')$, $\sigma \mapsto \sigma^{-1}$ is the identity on $MPR_{lsd}$.

[13] As a matter of fact the closed permutations that play a big role in [1] are maximal elements in their descent class (under the weak order), see also section 16 below.

In the case of the left weak order this works out as follows. There is an edge from $\sigma$ down to $\tau$ iff there are an $s_i, s_j$ with $i < j$ and $s_i = s_j + 1$ and $\tau$ is obtained from $\sigma$ by switching $s_i$ and $s_j$ (and leaving everything else as it was).

The Hasse diagrams for the left weak orders on $S_2$, $S_3$, and $S_4$, respectively look as follows

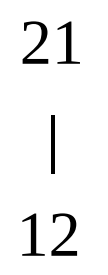

$S_2$

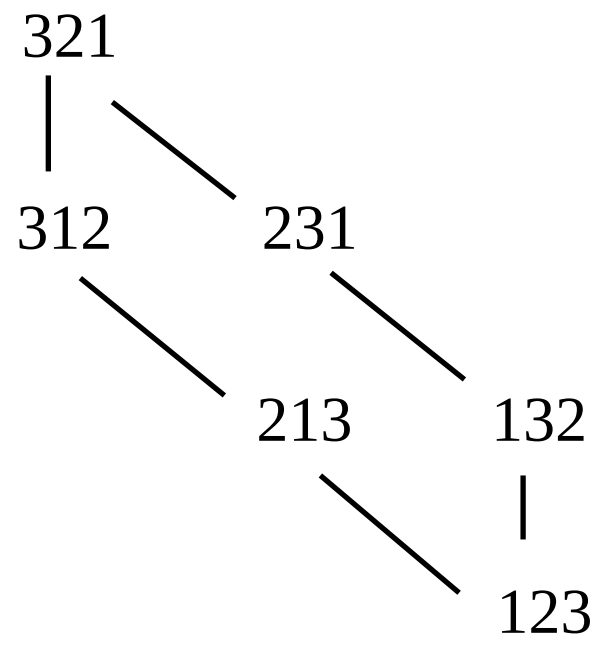

$S_3$

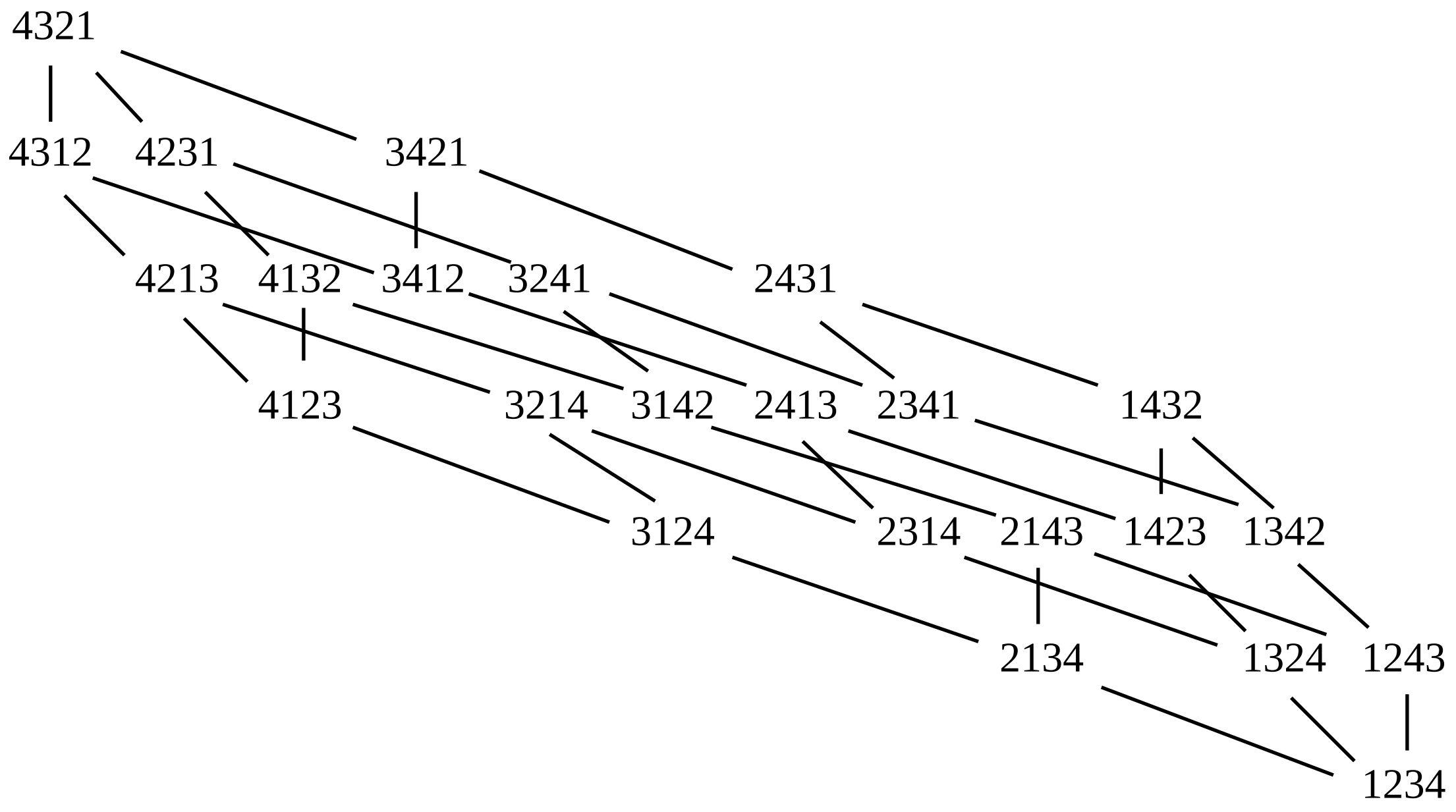

$S_4$

This is not the usual way of depicting the left weak order. It has been done this way here to bring out the recursive structure of the weak order.

The diagram for $S_n$ consists of $n$ copies of that for $S_{n-1}$ arranged in a sort of staircase, each step consisting of those elements that start with a fixed $i$.. For each $i$ the $i$-th step is the Hasse diagram of permutations of $\{1,2,\cdots,n\} \setminus \{i\}$ with for $i = 2,\cdots,n-1$ certain edges left out viz those that correspond to switches of $i-1$ and $i+1$ (which are neighbours in $\{1,2,\cdots,n\} \setminus \{i\}$ but not neighbours in $\{1,2,\cdots,n\}$). Finally corresponding members of step $j$ and step $j-1$ are connected for $j = n, n-1,\cdots,2$.

Note that on each layer the elements occur in lexicographic order from left to right. As will be seen this is by no means the only compatibility between the weak order and lexicographic order. Note also that the lexicographic order refines the weak order in the sense that $\sigma >_{lwo} \tau \quad \Rightarrow \quad \sigma >_{lex} \tau$. All this is practically immediate from from the description of the covering relation (Hasse diagram) of the weak order given above.

15.2. *Theorem.* Let $D \subset \{1,\cdots m-1\}$ be a descent set. Then the subgraph of the Hasse diagram of the left weak order formed by the descent class of $D$, i.e. by the permutations in $S_m$ with descent set $D$, form a connected subgraph with a unique smallest element (with respect to the weak order). That element is lsd. There is also a unique largest element in that descent class (with respect to the weak order) and this element is lld (**l**exicographically **l**argest in its **d**escent class).

*Proof.* This is proved by induction on length, i.e. by induction on $m$. Let $\sigma = [a_1,\cdots,a_m]$ be a permutation word. There are two cases to consider.

Case A: $1 \notin D$. Suppose that $a_1 \neq 1$. Because $a_2 > a_1$, $a_1 - 1 \in \{a_3,\cdots,a_m\}$. So switching $a_1$ and $a_1 - 1$ a permutation is obtained that is smaller in the weak order, which has the same descent set and which has smaller first element. Repeating this one finally obtains a permutation that is weak order smaller, that has the same descent set, and which has first element 1, say $\tau = [1,b_2,\cdots.b_m]$. Now by induction $[b_2,\cdots,b_m]$ can be changed by weak order decreasing switchings to an lsd permutation word (on the alphabeth $\{2,\cdots,m\}$ with the same descent set throughout. Then the same holds for $\tau$, because whatever is done to $[b_2,\cdots,b_m]$ never will 1 become a member of the descent set of the new permutation.

Case B. $\{1,\cdots,i-1\} \subset D$, $i \notin D$. (NB, if there is no $i \notin D$, $\sigma = [m,m-1,\cdots,2,1]$ and there is nothing to prove.) So $\sigma$ is of the form

$$\sigma = [a_1 > a_2 > \cdots > a_i < a_{i+1}, a_{i+2},\cdots,a_m]$$

So, if $a_i > 1$, it follows that $a_i - 1$ is in $\{a_{i+2},\cdots,a_m\}$. Thus $a_i$ and $a_i - 1$ can be interchanged to obtain a weak order smaller permutation with the same descent set and smaller number at place $i$. Continuing $\sigma$ can be brought in the form

$$\sigma = [a_1 > a_2 > \cdots > 1 < a_{i+1}, a_{i+2},\cdots,a_m]$$

by weak order decreasing transformation and with the same descent set throughout. With a subsidiary induction it can now be assumed that $\sigma$ has been brought into the form

$$\sigma = [a_1 > a_2 > \cdots > a_k > r > r-1 > \cdots > 1 < a_{i+1}, a_{i+2},\cdots,a_m]$$

Then $a_k - 1$ is not in $\{a_1,\cdots,a_{k-1}\}$ and also not in $\{1,\cdots,r\}$ unless $a_k = r+1$. Thus if $a_k \neq r+1$, switching $a_k$ and $a_k - 1$ produces a weak order smaller permutation with the same descent set and smaller element at spot $k$. Thus $\sigma$ can be transformed into

$$\sigma = [i,i-1,\cdots,2,1,a_{i+1},\cdots a_m] \tag{5.3}$$

Now by induction $[a_{i+1},\cdots,a_m]$ can be transformed by descent set preserving and weak order decreasing switchings into an lsd word $\tau$ on the alphabeth $\{i+1,\cdots,m\}$ and then

$$[i,i-1,\cdots,2,1] * \tau$$

has the same descent set as $\sigma$ and is also lsd. This proves the first two statements. The two statements concerning largest elements follow because the mapping

$$\sigma = [a_1, \cdots, a_m] \mapsto [m+1-a_1, \cdots, m+1-a_m] = \tau$$

reverses the weak order (and also the lexicographic order) while $\mathrm{desc}(\tau) = \{1, \cdots m-1\} \setminus \mathrm{desc}(\sigma)$.

**16. Global ascent and descent sets.**
An ascent of a permutation word $\sigma = [a_1, \cdots, a_m]$ is a $p \in \{1,2,\cdots,m-1\}$ for which $a_p < a_{p+1}$. The set of ascents of a permutation is denoted $\mathrm{asc}(\sigma)$. An ascent $p \in \{1,2,\cdots,m-1\}$ is global iff

$$a_i < a_j \quad \text{for all } i \le p,\ j \ge p+1 \tag{16.1}$$

The set of global ascents of a permutation is denoted $\mathrm{gasc}(\sigma)$. Obviously $\mathrm{gasc}(\sigma) \subset \mathrm{asc}(\sigma)$, but they can be equal. Similarly a descent $p \in \{1,2,\cdots,m-1\}$ is a global descent iff

$$a_i > a_j \quad \text{for all } i \le p,\ j \ge p+1 \tag{16.2}$$

As far as I know this notion first appeared in [1]. The permutations for which the global descent set is equal to the descent set are called closed permutations in [1].

16.3. *Theorem.* A permutation $\sigma = [a_1, \cdots, a_m]$ is lsd iff $\mathrm{gasc}(\sigma) = \mathrm{asc}(\sigma)$. A permutation is closed iff it is lld.

*Proof.* From the description of lsd permutations in section 13 it is immediate that lsd permutations have their global ascents sets equal to their ascent sets. Inversely assume that each ascent of $\sigma = [a_1, \cdots, a_m]$ is global. Let $i$ be the first ascent, so that the start of the permutation looks like

$$a_1 > a_2 > \cdots > a_{i-1} > a_i < a_{i+1}$$

Because the ascent is global the element $1$ cannot occur to the right of $a_i$ and it follows that the start of the permutation must look like

$$a_1 > a_2 > \cdots > a_{i-1} > 1 < a_{i+1}$$

If $i = 1$ induction finishes the proof. If $i > 1$, the element $2$ must occur somewhere and it cannot be to the right of $a_i = 1$ (because $a_{i-1} \ge 2$). Thus the start of the permutation must look like

$$a_1 > a_2 > \cdots > 2 > 1 < a_{i+1}$$

Continuing one finds that the start of the permutation must be of the form

$$i > i-1 > \cdots > 2 > 1 < a_{i+1}$$

and induction finishes the proof of the first statement of the theorem. The second statement follows by observing that if

$$\sigma = [a_1, \cdots, a_m] \quad \text{and} \quad \tau = [m+1-a_1, \cdots m+1-a_m]$$

$\mathrm{asc}(\sigma) = \mathrm{desc}(\tau)$, $\mathrm{gasc}(\sigma) = \mathrm{gdesc}(\tau)$ and that this transformation reverses lexicographical order.

Let $MPR_{lld}$ be the subgroup of $MPR$ spanned by the lld permutations. Obviously if $u*v$ is a cut of a closed permutation, i.e. one for which the descent set and global descent set are equal the same is true of $u$ and $v$ and their standardizations. So these are also lld by theorem 16.3 and it follows that $MPR_{lld}$ is a subcoalgebra of $(MPR,\mu)$. There results another coalgebra section of $MPR \xrightarrow{\pi} QSymm$ just as in the case of theorem 14.4.

Another way to obtain this one is to use theorem 14.4 in combination with theorem 4.1 from [11] which says that $\sigma = [a_1,\cdots,a_m] \mapsto \tau = [m+1-a_1,\cdots,m+1-a_m]$ is an automorphism of coalgebras of $(MPR,\mu)$.

**17. The incisive cut coalgebra.**

Consider again the free Abelian group with as basis all words over the natural numbers. This time consider the following (possible) comultiplication on a basis word $\alpha = [a_1,a_2,\cdots a_m]$

$$\mu(\alpha) = \sum_{i=0}^{m}[a_1,\cdots,a_i]\otimes[a_{i+1},\cdots a_m] + \sum_{j=1}^{m}\sum_{\substack{b_j+c_j=a_j \\ b_j,\, c_j>0}}[a_1,\cdots,a_{j-1},b_j]\otimes[c_j,a_{j+1},\cdots,a_m] \qquad (17.1)$$

The first sum consists of the normal cuts as in *QSymm*; in the second sum each $a_j$ that is greater than 1 is split in all possible ways into two parts, the incisive cuts. It is an easy exercise to check that this is coassociative and that with the same counit as *QSymm*, i.e. $\varepsilon$ takes the value 1 on the empty word and the value zero on all words of length greater than zero, this becomes a graded connected coalgebra. It will be denoted *ICC*.

17.2. *Theorem.* Consider the morphisms of Abelian groups defined by

$$(MPR,\mu) \longrightarrow ICC, \quad \sigma \mapsto \text{comp}(\text{desc}(\sigma) \qquad (17.3)$$

where comp(desc($\sigma$) is the composition associated to the descent set of the permutation $\sigma$. Then (17.3) is a morphism of coalgebras.

To start with, here is a reassuring remark that this may work out just right. The number of terms on the right hand side of (17.1) is obviously

$$m+1+\sum_{j=1}^{m}a_j - 1 = \text{wt}(\alpha)+1$$

and the number of terms in $\mu(\sigma)$ is $\text{lg}(\sigma)+1$, while the weight of comp(desc($\sigma$) is $\text{lg}(\sigma)$; so things fit.

*Proof of theorem* 17.1. Let $\text{lg}(\sigma) = n$ (i.e. $\sigma$ is a permutation on $n$ letters), and let its descent set be

$$D = \{d_1,\cdots,d_{m-1}\} \subset \{1,2,\cdots,n-1\}$$

so that $\alpha = $ comp(desc($\sigma$) has length $m$ and weight $n$. Recall that

$$\alpha = [a_1,\cdots,a_m], \quad a_i = d_i - d_{i-1}, \quad \text{with } d_0 = 0,\ d_m = n$$

Now let $\sigma = u*v$ be a cut with $u$ of length $r$ and $v$ of length $n-r$. The corresponding term in $\mu(\sigma)$ is $\text{st}(u)\otimes\text{st}(v)$. The question is now what are the compositions associated to the descent sets of $\text{st}(u)$ and $\text{st}(v)$ which are the same as those of $u$ and $v$. Let $j$ be the unique

index such that

$$d_j \le r < d_{j+1}, \quad j \in \{0,1,\cdots,m-1\}$$

There are two cases.

If $d_j = r$, the descent set of $u$ is equal to $\{d_1,\cdots,d_{j-1}\} \subset \{1,2,\cdots,r-1\}$, so that its composition is $[a_1,\cdots,a_j]$ while the descent set of $v$ is $\{d_{j+1}-r,\cdots,d_{m-1}-r\} \subset \{1,2,\cdots,n-r-1\}$ with associated composition $[a_{j+1},\cdots,a_m]$. So in this case there results a normal cut, a term from the first sum in (17.1).

If $d_j < r$, let $b_{j+1} = r - d_j$, $c_{j+1} = d_{j+1} - r$. Note that $b_{j+1} + c_{j+1} = a_{j+1}$. In this case the descent set of $u$ is equal to $\{d_1,\cdots,d_j\} \subset \{1,2,\cdots,r-1\}$ and that of $v$ is $\{d_{j+1}-r,\cdots,d_{m-1}-r\} \subset \{1,2,\cdots,n-r-1\}$. The associated compositions are $[a_1,\cdots,a_j,b_{j+1}]$ and $[c_{j+1},a_{j+2},\cdots,a_m]$, and so in this case there results an incisive cut, a term from the second sum in (17.1). Also all terms arise. This proves the theorem.

Consider the graded dual Abelian group $ICC^{gr*}$ of $ICC$ and let the dual basis be $\{R'_\alpha\}$ indexed by all words over the natural numbers. The coalgebra structure on $ICC$ induces an algebra structure on $ICC^{gr*}$ that on the basis $\{R'_\alpha\}$ is (obviously) given by

$$R'_{[a_1,\cdots,a_m]}R'_{[b_1,\cdots,b_n]} = R'_{[a_1,\cdots,a_m,b_1,\cdots,b_n]} + R'_{[a_1,\cdots,a_{m-1},a_m+b_1,b_2,\cdots,b_n]} \tag{17.4}$$

This looks familiar. Indeed, it is the multiplication of the ribbon Schur functions in *NSymm*, see (12.13). Thus $ICC^*$ is isomorphic as an algebra to *NSymm*, and hence $ICC$ is isomorphic as a coalgebra to *QSymm*. The latter isomorphism is given by $\alpha \mapsto F_\alpha$ because the $R_\alpha$ and the $F_\alpha$ form dual bases, see section 12.

This sheds, I feel, some extra light on the naturalness of the projection

$$\pi: \ MPR(m,\mu) \longrightarrow QSymm, \quad \sigma \mapsto F_{\text{comp}(\text{desc}(\sigma))}$$

**18. Coda.**

The study of such 'word Hopf algebras' as defined in this paper has just started (if that).

In particular many of them have large families of (Hopf) algebra endomorphisms (which is not, I believe, the case for *MPR*). These beg to be studied. Having Hopf algebras with lots of endomorphisms was and is one of the motivations for introducing word Hopf algebras.

All of the word Hopf algebras above are based (loosely speaking) on endomorphism recipes for *Shuffle* (and *LieHopf*). What happens when endomorphism recipes are considered for other (connected graded) Hopf algebras such as *NSymm*, *QSymm*, *MPR*, and even the word Hopf algebras themselves?

**References.**

1. M. Aguiar and F. Sottile, *Structure of the Malvenuto-Reutenauer Hopf algebra of permutations*, Texas A&M University, 2002, to appear Adv. Math.

2. R. J. Blattner, *Crossed products and inner actions of Hopf algebras*, Trans. Amer. Math. Soc. **298** (1986), 671 - 711.

3. Y. Doi and M. Takeuchi, *Cleft comodule algebras for a bialgebra*, Comm. Algebra, **14** (1986), 801-818.

4. G. Duchamp, F. Hivert and J.-Y. Thibon, *Noncommutative symmetric functions VI: free quasi-symmetric functions and related algebras*, Université de Marne-la-Vallée, 2001.

5. I. M. Gelfand, D. Krob, A. Lascoux, B. Leclerc, V. S. Retakh and Jean-Yves-Thibon, *Noncommutative symmetric functions*, Adv. Math. **112** (1995), 218-348.

6. M. Hazewinkel, *The primitives of the Hopf algebra of noncommutative symmetric functions over the integers*, CWI, 2001, submitted Ann. Math.

7. M. Hazewinkel, *Symmetric functions, noncommutative symmetric functions, and quasisymmetric functions*, Acta Appl. Math. **75** (2003), 55-83.

8. M. Hazewinkel, *Symmetric functions, noncommutative symmetric functions, and quasisymmetric functions II*, Acta Appl. Math. (2004, to appear).

9. M. Hazewinkel, N. Gubareni and V. Kirichenko, *Algebras, rings, and modules. Volume 2*, KAP, to appear.

10. C. Malvenuto and C. Reutenauer, *Duality between quasi-symmetric functions and the Solomon descent algebra*, J. of Algebra **177** (1994), 967-982.

11. S. Poirier and C. Reutenauer, *Algèbres de Hopf de tableaux*, Ann. Sci. Math. Québec **19** (1995), 79-90.

12. C. Reutenauer, *Free Lie algebras*, Oxford University Press, 1993.

13. C. Schensted, *Longest increasing and decreasing subsequences*, Canadian J. Math. **13** (1961), 179-191.

14. S. Shnider and S. Sternberg, *Quantum groups. From coalgebras to Drinfeld algebras*, International Press, 1994.

15. L. Solomon, *A Mackey formula in the group ring of a Coxeter group*, J. of Algebra **41** (1976), 255-268.